\documentclass[12pt]{amsart}
\usepackage{amssymb,amscd}
\usepackage{verbatim}

\usepackage{amsmath,amssymb,graphicx,mathrsfs}   
\usepackage[colorlinks=true,allcolors = blue]{hyperref} 

\let\frak\mathfrak
\let\Bbb\mathbb

\def\>{\relax\ifmmode\mskip.666667\thinmuskip\relax\else\kern.111111em\fi}
\def\<{\relax\ifmmode\mskip-.333333\thinmuskip\relax\else\kern-.0555556em\fi}
\def\vsk#1>{\vskip#1\baselineskip}
\def\vv#1>{\vadjust{\vsk#1>}\ignorespaces}
\def\vvn#1>{\vadjust{\nobreak\vsk#1>\nobreak}\ignorespaces}
\def\vvgood{\vadjust{\penalty-500}} \let\alb\allowbreak
\def\fratop{\genfrac{}{}{0pt}1}
\def\satop#1#2{\fratop{\scriptstyle#1}{\scriptstyle#2}}
\let\dsize\displaystyle \let\tsize\textstyle \let\ssize\scriptstyle
\let\sssize\scriptscriptstyle \def\tfrac{\textstyle\frac}
\let\phan\phantom \let\vp\vphantom \let\hp\hphantom

\def\sskip{\par\vsk.2>}
\let\Medskip\medskip
\def\medskip{\par\Medskip}
\let\Bigskip\bigskip
\def\bigskip{\par\Bigskip}

\let\Maketitle\maketitle
\def\maketitle{\Maketitle\thispagestyle{empty}\let\maketitle\empty}

\newtheorem{thm}{Theorem}[section]
\newtheorem{cor}[thm]{Corollary}
\newtheorem{lem}[thm]{Lemma}

\newtheorem{exmp}[thm]{Example}

\theoremstyle{definition}                                  
\numberwithin{equation}{section}

\theoremstyle{definition}
\newtheorem*{rem}{Remark}

\let\mc\mathcal
\let\nc\newcommand

\def\flati{\def\=##1{\rlap##1\hphantom b)}}

\let\al\alpha
\let\bt\beta
\let\dl\delta
\let\Dl\Delta
\let\eps\varepsilon
\let\gm\gamma
\let\Gm\Gamma
\let\ka\kappa
\let\la\lambda
\let\La\Lambda
\let\pho\phi
\let\phi\varphi
\let\si\sigma
\let\Si\Sigma
\let\Sig\varSigma
\let\Tht\Theta
\let\tht\theta
\let\thi\vartheta
\let\Ups\Upsilon
\let\om\omega
\let\Om\Omega

\let\der\partial
\let\Hat\widehat
\let\minus\setminus
\let\ox\otimes
\let\Tilde\widetilde
\let\bra\langle
\let\ket\rangle

\let\ge\geqslant
\let\geq\geqslant
\let\le\leqslant
\let\leq\leqslant

\let\on\operatorname
\let\bi\bibitem
\let\bs\boldsymbol
\let\xlto\xleftarrow
\let\xto\xrightarrow

\def\C{{\mathbb C}}
\def\Z{{\mathbb Z}}
\def\R{{\mathbb R}}
\def\PP{{\mathbb P}}

\def\Ac{{\mc A}}
\def\B{{\mc B}}
\def\F{{\mathbb F}}   

\def\Hc{{\mc H}}
\def\Jc{{\mc J}}
\def\Oc{{\mc O}}
\def\Sc{{\mc S}}
\def\V{{\mc V}}

\def\+#1{^{\{#1\}}}
\def\lsym#1{#1\alb\dots\relax#1\alb}
\def\lc{\lsym,}
\def\lox{\lsym\ox}

\def\diag{\on{diag}}
\def\End{\on{End}}
\def\Gr{\on{Gr}}
\def\Hom{\on{Hom}}
\def\id{\on{id}}
\def\Res{\on{Res}}
\def\qdet{\on{qdet}}
\def\rdet{\on{rdet}}
\def\rank{\on{rank}}
\def\sign{\on{sign}}
\def\tbigoplus{\mathop{\textstyle{\bigoplus}}\limits}
\def\tbigcap{\mathrel{\textstyle{\bigcap}}}
\def\tr{\on{tr}}
\def\Wr{\on{Wr}}

\def\ii{i,\<\>i}
\def\ij{i,\<\>j}
\def\ik{i,\<\>k}
\def\il{i,\<\>l}
\def\ji{j,\<\>i}
\def\jj{j,\<\>j}
\def\jk{j,\<\>k}
\def\kj{k,\<\>j}
\def\kl{k,\<\>l}
\def\lk{l,\<\>k}
\def\II{I\<\<,\<\>I}
\def\IJ{I\<\<,J}
\def\ioi{i+1,\<\>i}
\def\pp{p,\<\>p}
\def\ppo{p,\<\>p+1}
\def\pop{p+1,\<\>p}
\def\pci{p,\<\>i}
\def\pcj{p,\<\>j}
\def\poi{p+1,\<\>i}
\def\poj{p+1,\<\>j}

\def\gl{\mathfrak{gl}}
\def\gln{\mathfrak{gl}_N}
\def\sln{\mathfrak{sl}_N}
\def\slnn{\mathfrak{sl}_{\>n}}

\def\glnt{\gln[t]}
\def\Ugln{U(\gln)}
\def\Uglnt{U(\glnt)}
\def\Yn{Y\<(\gln)}

\def\slnt{\sln[t]}
\def\Usln{U(\sln)}
\def\Uslnt{U(\slnt)}

\def\beq{\begin{equation}}
\def\eeq{\end{equation}}
\def\be{\begin{equation*}}
\def\ee{\end{equation*}}

\nc{\bea}{\begin{eqnarray*}}
\nc{\eea}{\end{eqnarray*}}
\nc{\bean}{\begin{eqnarray}}
\nc{\eean}{\end{eqnarray}}
\nc{\Ref}[1]{{\rm(\ref{#1})}}

\def\g{{\mathfrak g}}
\def\h{{\mathfrak h}}
\def\n{{\mathfrak n}}

\let\ga\gamma
\let\Ga\Gamma

\nc{\Il}{{\mc I_{\bs\la}}}
\nc{\bla}{{\bs\la}}
\nc{\Fla}{\F_\bla}
\nc{\tfl}{{T^*\Fla}}
\nc{\GL}{{GL_n(\C)}}
\nc{\GLC}{{GL_n(\C)\times\C^*}}

\let\aal\al 
\def\mub{{\bs\mu}}
\def\Dc{\check D}
\def\Di{{\tfrac1D}}
\def\Dci{{\tfrac1\Dc}}
\def\DV{\Di\V^-}
\def\DL{\DV_\bla}
\def\Vl{\V^+_\bla}
\def\DVe{\Dci\V^=}
\def\DLe{\DVe_\bla}

\def\xxx{x_1\lc x_n}
\def\yyy{y_1\lc y_n}
\def\zzz{z_1\lc z_n}
\def\Cx{\C[\xxx]}
\def\Czh{\C[\zzz,h]}
\def\Vz{V\<\ox\Czh}
\def\ty{\Tilde Y\<(\gln)}
\def\tb{\Tilde\B}
\def\IMA{{I^{\<\>\max}}}
\def\IMI{{I^{\<\>\min}}}
\def\CZH{\C[\zzz]^{\>S_n}\!\ox\C[h]}
\def\Sla{S_{\la_1}\!\lsym\times S_{\la_N}}
\def\SlN{S_{\la_N}\!\lsym\times S_{\la_1}}
\def\Czhl{\C[\zzz]^{\>\Sla}\!\ox\C[h]}
\def\Czghl{\C[\<\>\zb\<\>;\Gmm\>]^{\>S_n\times\Sla}\!\ox\C[h]}
\def\Czs{\C[\zzz]^{\>S_n}}

\def\Ct{\Tilde C}
\def\fc{\check f}
\def\Hg{{\mathfrak H}}
\def\Ic{\check I}
\def\Ih{\hat I}
\def\It{\tilde I}
\def\Jt{\tilde J}
\def\Pit{\Tilde\Pi}
\def\Qc{\check Q}
\let\sd s 
\def\sh{\hat s}
\def\st{\tilde s}
\def\sih{\hat\si}
\def\sit{\tilde\si}
\def\nablat{\Tilde\nabla}
\def\Ut{\Tilde U}
\def\Xt{\Tilde X}
\def\Dh{\Hat D}
\def\Vh{\Hat V}
\def\Wh{\Hat W}

\def\ib{\bs i}
\def\jb{\bs j}
\def\kb{\bs k}
\def\iib{\ib,\<\>\ib}
\def\ijb{\ib,\<\>\jb}
\def\ikb{\ib,\<\>\kb}
\def\jib{\jb,\<\>\ib}
\def\jkb{\jb,\<\>\kb}
\def\kjb{\kb,\<\>\jb}
\def\zb{\bs z}
\def\zzi{z_1\lc z_i,z_{i+1}\lc z_n}
\def\zzii{z_1\lc z_{i+1},z_i\lc z_n}
\def\zzzn{z_n\lc z_1}
\def\xxi{x_1\lc x_i,x_{i+1}\lc x_n}
\def\xxii{x_1\lc x_{i+1},x_i\lc x_n}
\def\ip{\<\>i\>\prime}
\def\ipi{\>\prime\<\>i}
\def\jp{\prime j}
\def\iset{\{\<\>i\<\>\}}
\def\jset{\{\<\>j\<\>\}}
\def\IMIp{{I^{\min,i\prime}}}
\let\Gmm\Gamma
\def\top{\on{top}}
\def\phoo{\pho_{\>0}}
\def\Vy#1{V^{\bra#1\ket}}
\def\Vyh#1{\Vh^{\bra#1\ket}}
\def\Hh{\Hg_n}
\def\bmo{\nabla^{\>\sssize\mathrm{BMO}}}

\def\Bin{B^{\<\>\infty}}
\def\Bci{\B^{\<\>\infty}}
\def\tbi{\tb^{\<\>\infty}}
\def\Xin{X^\infty}

\def\Bb{\rlap{$\<\>\bar{\phan\B}$}\rlap{$\>\bar{\phan\B}$}\B}
\def\Hb{\rlap{$\bar{\,\,\phan H}$}\bar H}
\def\Hbc{\rlap{$\<\>\bar{\phan\Hc}$}\rlap{$\,\bar{\phan\Hc}$}\Hc}
\def\Xb{\rlap{$\bar{\;\>\phan X}$}\rlap{$\bar{\<\phan X}$}X}

\def\Bk{B^{\<\>\kk}}
\def\Bck{\B^{\<\>\kk}}
\def\Bbk{\Bb^{\<\>\kk}}
\def\Hck{\Hc^{\<\>\kk}}
\def\Hbck{\Hbc^{\<\>\kk}}
\def\mukp{\mu^{\kk+}}
\def\muke{\mu^{\kk=}}
\def\mukm{\mu^{\kk-}}
\def\mukpm{\mu^{\kk\pm}}
\def\nukp{\nu^{\<\>\kk+}}
\def\nuke{\nu^{\<\>\kk=}}
\def\nukm{\nu^{\<\>\kk-}}
\def\nukpm{\nu^{\<\>\kk\pm}}
\def\rhkm{\rho^{\<\>\kk-}}
\def\rhkpm{\rho^{\<\>\kk\pm}}
\def\tbk{\tb^{\<\>\kk}}
\def\tbkp{\tb^{\<\>\kk+}}
\def\tbkm{\tb^{\<\>\kk-}}
\def\tbkpm{\tb^{\<\>\kk\pm}}
\def\tbbp{\Bb^{\<\>\kk+}}
\def\tbbm{\Bb^{\<\>\kk-}}
\def\Sk{S^{\>\kk}}
\def\Uk{U^\kk}
\def\Upk{U'^{\<\>\kk}}
\def\Wk{W^\kk}
\def\Whk{\Wh^\kk}
\def\Xk{X^\kk}
\def\Xkp{X^{\kk+}}
\def\Xkm{X^{\kk-}}
\def\Xkpm{X^{\kk\pm}}
\def\Xtk{\Xt^\kk}
\def\Xbk{\Xb^{\<\>\kk}}
\def\Yk{Y^{\<\>\kk}}
\def\ddk_#1{\kk_{#1}\<\>\frac\der{\der\<\>\kk_{#1}}}
\def\zno{\big/\bra\>z_1\<\lsym+z_n\<=0\,\ket}
\def\vone{v_1\<\lox v_1}
\def\Hplus{\bigoplus_\bla H_\bla}
\def\Pone{\mathbb P^{\<\>1}}
\def\FFF{\mathbb{F}}

\def\bul{\mathbin{\raise.2ex\hbox{$\sssize\bullet$}}}
\def\intt{\mathchoice
{\mathop{\raise.2ex\rlap{$\,\,\ssize\backslash$}{\intop}}\nolimits}
{\mathop{\raise.3ex\rlap{$\,\sssize\backslash$}{\intop}}\nolimits}
{\mathop{\raise.1ex\rlap{$\sssize\>\backslash$}{\intop}}\nolimits}
{\mathop{\rlap{$\sssize\<\>\backslash$}{\intop}}\nolimits}}

\def\nablas{\nabla^{\<\>\star}}
\def\nablab{\nabla^{\<\>\bul}}
\def\chib{\chi^{\<\>\bul}}

\def\zla{\zeta_{\>\bla}^{\vp:}}
\let\gak\gamma 
\let\kk q 
\let\kp\kappa 
\let\cc c
\def\kkk{\kk_1\lc\kk_N}
\def\kkn{\kk_1\lc\kk_n}
\let\Ko K
\def\Kh{\Hat\Ko}
\def\zzip{z_1\lc z_{i-1},z_i\<-\kp,z_{i+1}\lc z_n}
\def\prodl{\prod^{\longleftarrow}}
\def\prodr{\prod^{\longrightarrow}}

\def\GZ/{Gelfand-Zetlin}
\def\KZ/{{\slshape KZ\/}}
\def\qKZ/{{\slshape qKZ\/}}
\def\XXX/{{\slshape XXX\/}}

\def\zz{{\bs z}}
\def\qq{{\bs q}}
\def\TT{{\bs t}}
\def\glN{{\frak{gl}_N}}
\def\ts{\bs t}
\def\Sym{\on{Sym}}
\def\Cts{\C[\<\>\ts\<\>]^{\>S_{\la^{(1)}}\<\<\lsym\times S_{\la^{(N-1)}}}}
\def\ss{{\bs s}}
\def\BL{{\tilde{\mc B}^q(H_\bla)}}
\nc{\A}{{\mc C}}
\def\glnn{{\frak{gl}_n}}
\def\uu{{\bs u}}
\def\VV{{\bs v}}
\def\GG{{\bs \Ga}}
\def\II{{\mc I}}
\def\XX{{\mc X}}
\def\a{{\frak{a}}}
\def\CC{{\frak{C}}}
\def\St{{\on{Stab}}}
\def\xx{{\bs x}}
\def\Czh{{(\C^N)^{\otimes n}\otimes\C(\zz;h)}}

\def\WW{{\check{W}}}
\def\FF{{\mathbb F}}
\def\Sing{{\on{Sing}}}
\def\sll{{\frak{sl}}}

\def\slt{{\frak{sl}_2}}
\def\Ant{{\on{Ant}}}
\def\Ik{{\mc I_k}}
\def\Q{{\mathbb Q}}
\def\K{{\mathbb K}}

\begin{document}

\hrule width0pt
\vsk->
\title[Hypergeometric Integrals Modulo $p$]{Hypergeometric Integrals Modulo $p$ \\  and Hasse--Witt Matrices}

\author[Alexey Slinkin and Alexander Varchenko]
{Alexey Slinkin$^{\diamond}$ and Alexander Varchenko$^{\star}$}

\maketitle

\begin{center}
{\it $^{\diamond, \star}\<$Department of Mathematics, University
of North Carolina at Chapel Hill\\ Chapel Hill, NC 27599-3250, USA\/}

\vsk.5>
{\it $^{\star}\<$Faculty of Mathematics and Mechanics, Lomonosov Moscow State
University\\ Leninskiye Gory 1, 119991 Moscow GSP-1, Russia\/}

\end{center}

\vsk>
{\leftskip3pc \rightskip\leftskip \parindent0pt \Small
{\it Key words\/}:  KZ equations, hypergeometric integrals, Hasse-Witt matrix, reduction to
\\
\phantom{aaaaaaaaaa}
characteristic $p$

\vsk.6>
{\it 2010 Mathematics Subject Classification\/}: 13A35 (33C60, 32G20) 
\par}

{\let\thefootnote\relax
\footnotetext{\vsk-.8>\noindent
$^\diamond\<$
{\it E\>-mail}: slinalex@live.unc.edu
\\
$^\star\<$
{\it E\>-mail}:
anv@email.unc.edu\,, supported in part by NSF grant DMS-1665239}}

\begin{abstract} 
We consider the KZ differential equations over $\C$
in the case, when the hypergeometric solutions are 
one-dimensional integrals. We also consider
the same differential equations over a finite field $\F_p$.
We study  the space of polynomial solutions of these differential equations over $\F_p$, 
constructed in a previous work by
V.\,Schechtman and the second author.
Using Hasse-Witt matrices we identify the space of these polynomial solutions over $\F_p$ 
with the space dual to a certain subspace of regular differentials on an associated curve. We also 
relate these polynomial solutions over $\F_p$ and the hypergeometric solutions over $\C$.

\end{abstract}

{\small\tableofcontents\par}

\setcounter{footnote}{0}
\renewcommand{\thefootnote}{\arabic{footnote}}

\section{Introduction}

The KZ equations were discovered by  Vadim Knizhnik and Alexander Zamolodchikov \cite{KZ} to describe the differential equations for conformal blocks on sphere in the Wess-Zumino-Witten model of conformal field theory.
The hypergeometric solutions of the KZ equations were constructed more than  30 years ago,
see \cite{SV1, SV2}. The polynomial
solutions of the KZ equations over the finite field $\F_p$  
with a prime number $p$ of elements were constructed 
recently in \cite{SV3}.  We call these  solutions over $\F_p$ the {\it arithmetic solutions}.
The general problem is to find the dimension of the space of arithmetic solutions and to 
understand relations between the hypergeometric solutions 
of the KZ equations over $\C$ and the arithmetic solutions over $\F_p$.

\vsk.2>
In this paper we consider an example of the KZ differential equations, whose 
hypergeometric solutions over  $\C$ are $n$-vectors of  the integrals
\bean
\label{intr in}
I^{(\ga)}(z_1,\dots,z_n) = 
\left(\int_{\ga} \frac 1{x-z_1} \frac {dx}y,\ \dots, \ \int_{\ga} \frac 1{x-z_n} \frac {dx}y\right),
\eean
where
\bean
\label{intr y}
y^q = (x-z_1)\dots(x-z_n)\,
\eean
and $\ga$ is a suitable 1-cycle.
It is well known that the space of such $n$-vectors  is $n-1$-dimensional.
We consider the same  differential KZ equations over the field $\F_p$
under the assumption that $q$ is also a prime number and $p>q$, $p> n$, $n=kq+1$ for some positive integer $k$.
We show  that the dimension of the space of arithmetic solutions equals only a fraction of $n$.
Namely, let $a_1$ be the unique positive integer such that $1\leq a_1\leq q-1$ and
$a_1p\equiv 1$ (mod $q$). This $a_1$ is the inverse of $p$ modulo $q$. It turns out that 
the dimension of the space of arithmetic solutions equals 
$a_1k$. More precisely the dimension of the space of arithmetic solutions can be defined as follows.
Consider the curve $X$ defined by the affine equation \Ref{intr y}.
The cyclic group $\Z_q$ of $q$-th roots of unity acts on $X$ by multiplication on the coordinate $y$.
The space $\Om^1(X)$ of regular differentials on $X$ splits into eigenspaces of the $\Z_q$-action,
$\Omega^1 (X) \,=\, \bigoplus_{a=1}^{q-1}\, \Omega^1_a (X)\,,$ where
$ \Omega^1_a (X)$ consists of differentials of the form $u(x) dx/y^a$. We show that the dimension
of the space of arithmetic solutions equals the dimension of $\Om^1_{a_1}(X)$.
Moreover, we establish an isomorphism of the space of arithmetic solutions  and the space dual to 
$\Om^1_{a_1}(X)$. That isomorphism is constructed with the help of the map adjoint to the
corresponding Cartier map, and more precisely, with the help of the corresponding  Hasse-Witt matrix.
This is our first main result.

\vsk.2>
We also choose one  solution 
of the KZ equations over $\C$ and call it {\it distinguished}. We expand 
the distinguished solutions into the Taylor series at some point, reduce the coefficients
of the Taylor expansion modulo $p$ and present this reduced Taylor series as an infinite formal sum of 
arithmetic solutions,  with coefficients being matrix elements of  iterates of the associated Hasse-Witt matrix.
Moreover, this presentation allows one to recover a basis of arithmetic solutions in terms of the reduced Taylor expansion.
This statement is our second main result.  

Our comparison of the reduced Taylor expansion of a solution over $\C$
and arithmetic solutions over $\F_p$ is analogous to 
Y.I.\,Manin\rq{}s considerations of the elliptic integral in his classical paper \cite{Ma} in 1961,
 see also Section ``Manin\rq{}s Result: The Unity of Mathematics\rq{}\rq{}
 in the book \cite{Cl} by Clemens.

\vsk.2>
For $q=2$ the results of this paper have been obtained in \cite{V5}.

\vsk.2>

The paper is organized as follows. In Section \ref{sec KZ} we describe the differential
KZ equations considered in this paper. 
We construct the hypergeometric solutions of these equations over $\C$ and the arithmetic solutions over $\F_p$.
In Section \ref{sec mod M} we define the module of arithmetic solutions, 
show that it is free, and calculate the rank of the module.
 In Section \ref{sec fus} we describe the fusion procedure for modules of arithmetic solutions.  

In Section \ref{sec bin coef} we prove a generalization of the classical
Lucas theorem, which will allow us to reduce modulo $p$
the coefficient of the Taylor expansion of the distinguished solution.

In Section \ref{sec m=1} we discuss different bases in the 
module of arithmetic solutions. One of the bases naturally appears 
in the defining construction of arithmetic solutions and the
other is convenient to relate 
the arithmetic solutions to the Taylor expansion of the distinguished solution.

In Sections \ref{sec Cartier-Manin} and \ref{sec C mg} we study the Cartier map related to our 
arithmetic solutions and prove our first main result by
identifying the module of arithmetic solutions with the
space dual to $\Om^1_{a_1}(X)$, see Theorems \ref{thm C ar1} and \ref{thm C ar2}.
In Section \ref{sec HW y} we change variables in the Hasse-Witt matrix preparing it for application to the
study of the 
distinguished solution.

In Section \ref{sec 8} we compare the reduced Taylor series of the distinguished solution and arithmetic solutions,
see Theorem \ref{thm dec}.

\vsk.2>
The authors thank P.\,Etingof for useful discussions.

\section{KZ equations}
\label{sec KZ}

\subsection{Description of equations}   Let $\g$ be a simple Lie algebra  over the field $\C$,
$\Om \in\g^{\ox 2} $ the Casimir element corresponding to an invariant scalar product on $\g$, \
$V_{1}, \dots, V_{n}$ finite-dimensional irreducible $\g$-modules.

The system of KZ equations with parameter $\ka\in \C^\times$
on a tensor $\ox_{i=1}^nV_{i}$ valued function 
$I(z_1,\dots,z_n)$ is the system of the differential equations 
\bean
\label{KZg}
 \frac{\der I}{\der z_i} = 
\frac 1\ka \sum_{j\ne i}\frac{\Om^{(i,j)}}{z_i-z_j} I,\qquad i=1,\dots,n,
\eean
where $\Om^{(i,j)}$ is the Casimir element acting in the $i$-th and $j$-th factors,
see \cite{KZ, EFK}.  The KZ differential equations commute with the action of $\g$ on 
$\ox_{i=1}^nV_{i}$, in particular, they preserve the subspaces of singular vectors
of a given weight.

In \cite{SV1, SV2} the KZ equations restricted to the subspace of singular vectors of a
given weight were identified with 
a suitable Gauss-Manin differential equations and the corresponding solutions of the KZ equations were
presented as multidimensional hypergeometric integrals. 

\vsk.2>

Let $p$ be a prime number and $\F_p$ the field with $p$ elements. Let $\g^p$ be the same Lie algebra considered over 
 $\F_p$. 
Let $V_{1}^p, \dots, V_{n}^p$ be the $\g^p$-modules which are reductions modulo $p$ of
$V_{1}, \dots, V_{n}$, respectively.
If $\ka$ is an integer and $p$ large enough with respect to $\ka$, then one can look for solutions
$I(z_1,\dots,z_n)$ of the KZ equations in $\ox_{i=1}^nV_{i}^p\ox \F_p[z_1,\dots,z_n]$.
Such solutions were constructed in \cite{SV3}.

In this paper we address two questions: 
\begin{enumerate}
\item[A.]
What is the number of independent solutions constructed in \cite{SV3} for a given $\F_p$?
\item[B.]
How are those solutions related to the solutions over $\C$,  that 
are given by hypergeometric integrals?
\end{enumerate}

We answer these questions in an example in which the hypergeometric solutions are presented by one-dimensional
 integrals. The case of hyperelliptic integrals was considered in \cite{V5}.
 The object of our study is the following joint system of differential and algebraic 
equations. 

\subsubsection{Assumptions in this paper}
\label{sec ass}

{\it
We fix prime numbers $p, q$, $p>q$,  a positive integer $n$,
a vector $\La=(\La_1,\dots, \La_n) \in \Z^n_{>0}$, such that 
$\La_i<q$ for all $i=1,\dots,n$. }

\vsk.2>

For 
 $z=(z_1,\dots, z_{n})$, we study the
 column vectors  $I(z)=(I_1(z)$, \dots, $I_{n}(z))$ satisfying  the  system of
 differential and algebraic linear equations:
\bean
\label{KZ}
\phantom{aaa}
 \frac{\partial I}{\partial z_i} \ = \
   {\frac 1q} \sum_{j \ne i}
   \frac{\Omega_{ij}}{z_i - z_j}  I ,
\quad i = 1, \dots , n,
\qquad
\La_1I_1(z)+\dots+\La_nI_{n}(z)=0,
\eean
where
\bean
\label{Om_ij_reduced}
 \Omega_{ij} \ = \ \begin{pmatrix}
             & \vdots^i &  & \vdots^j &  \\
        {\scriptstyle i} \cdots & {-\La_j} & \cdots &
            \La_j   & \cdots \\
                   & \vdots &  & \vdots &   \\
        {\scriptstyle j} \cdots & \La_i & \cdots & {-\La_i}&
                 \cdots \\
                   & \vdots &  & \vdots &
                   \end{pmatrix} ,
\eean                    
and all other entries are zero.

\vsk.2>
In this paper this joint system of {\it differential and 
algebraic equations} will be called the {\it KZ differential equations}. 

\vsk.2>
We will  construct solutions of these KZ differential equations over  $\C$ and over  $\F_p$
and compare the properties of solutions.

\vsk.2>
\begin{rem}

The system of equations \Ref{KZ} is the system of the true KZ differential equations \Ref{KZg} with parameter $\ka=q$, 
associated with the Lie algebra $\sll_2$ and the subspace of singular vectors of weight $\sum_{i=1}^n\La_i-2$ of the tensor 
product
$V_{\La_1}\ox\dots\ox V_{\La_n}$, 
where $V_{\La_i}$ is the irreducible $\La_i+1$ dimensional  $\sll_2$-module, up to a gauge transformation, see 
this example in  \cite[Section 1.1]{V2}.

Notice also that the assumption $\La_i<q$ for $i=1,\dots,n$ appears, when one is interested in differential equations for 
$\sll_2$ conformal blocks with central charge $q-2$.
\end{rem}

\subsection{Solutions  over $\C$}

Consider the {\it master function}
\bean
\label{mast f}
\Phi(t,z_1,\dots,z_{n}) = \prod_{a=1}^{n}(t-z_a)^{-\La_a/q}
\eean
and  the ${n}$-vector  of hypergeometric  integrals
\bean
\label{Iga}
I^{(\ga)} (z)=(I_1(z),\dots,I_n(z)),
\eean
 where
\bean
\label{s}
I_j=\int  \Phi(t,z_1,\dots,z_{n})\, \frac {dt}{t-z_j},\qquad j=1,\dots,{n}\,.
\eean
The integrals $I_j$, $j=1,\dots,n$, are over an element $\ga$ of the first homology group
 of the algebraic curve with affine equation
\bea
y^q = (t-z_1)^{\La_1}\dots (t-z_{n})^{\La_n}\,.
\eea
Starting from such $\ga$,  chosen for given $\{z_1,\dots,z_{n}\}$, the vector $I^{(\ga)}(z)$ can be analytically continued as a multivalued holomorphic function of $z$ to the complement in $\C^n$ to the union of the
diagonal  hyperplanes $z_i=z_j$.

\begin{thm}
\label{thm1.1}

 The vector $I^{(\ga)}(z)$ satisfies  the KZ differential  equations  \Ref{KZ}.
\end{thm}

Theorem \ref{thm1.1} is a classical statement. 
Much more general algebraic and differential equations satisfied by analogous multidimensional hypergeometric integrals were considered in \cite{SV1, SV2}.  Theorem \ref{thm1.1} is discussed as an example in  \cite[Section 1.1]{V2}.

\begin{thm} [{\cite[Formula (1.3)]{V1}}]
\label{thm dim}

All solutions of  equations \Ref{KZ} have this form. 
Namely, the complex vector space of solutions of the form \Ref{Iga}-\Ref{s} is $n-1$-dimensional.

\end{thm}

This theorem follows from the determinant formula for multidimensional hypergeometric integrals  in    \cite{V1}, in particular,
from \cite[Formula (1.3)]{V1}.

\subsection{Solutions  over $\F_p$}

Polynomial solutions of system \Ref{KZ}, considered over the field $\F_p$, 
were constructed in \cite{SV3}.

\vsk.2>

For $i=1,\dots,n$, choose positive integers $M_i$  such that 
\bean
\label{M_i}
M_i\equiv -\frac{\La_i}{q}
\qquad 
(\on{mod} \,p)\,,
\eean
that is,  project $\La_i, q$ to $\FFF_p$, calculate $-\frac{\La_i}q$ in $\FFF_p$ and 
then choose positive integers $M_i$ satisfying these equations.  Denote $M=(M_1,\dots,M_n)$.
Consider the {\it master polynomial}
\bean
\label{mp red}
\Phi_{p}(t,z, M) := \prod_{i=1}^n(t-z_i)^{M_i},
\eean
and the Taylor expansion with respect to the variable $t$ of the vector of polynomials
\bea
P(t,z,M):=\Phi_{p}(t,z, M)\Big(\frac 1{t-z_1}, \dots,\frac 1{t-z_n}\Big)\,=\, \sum_i P^{i}(z, M) \,t^i,
\eea
where $P^{i}(z,M)$ are $n$-vectors of polynomials in $z_1,\dots,z_n$ with coefficients in $\FFF_p$.

\begin{thm}[{\cite[Theorem 1.2]{SV3}}] \label{thm Fp} 
For any positive integer $l$, the vector of polynomials
\\
 $P^{lp-1}(z,M)$ 
satisfies the KZ differential equations \Ref{KZ}.

\end{thm}

Theorem \ref{thm Fp} is a particular case of \cite[Theorem 2.4]{SV3}. Cf. Theorem \ref{thm Fp} in \cite{K}. See also \cite{V3, V4, V5}. 
\vsk.2>

The solutions  $P^{lp-1}(z,M)$ 
given by this construction will be called the {\it arithmetic solutions} of the KZ differential equations \Ref{KZ}.

\section{Module  of arithmetic solutions}
\label{sec mod M}

\subsection{Definition of the module}
\label{sec dM}

Denote $\F_p[z^p]:=\F_p[z_1^p,\dots,z_{n}^p]$. 
The set of all polynomial solutions of \Ref{KZ} with coefficients in $\F_p$
is a module over the ring
$\F_p[z^p]$ since equations \Ref{KZ} are linear and
$\frac{\der z_i^p}{\der z_j} =0$ in $\F_p[z]$ for all $i,j$.
\vsk.2>

The arithmetic solutions $P^{lp-1}(z,M)$  of equations \Ref{KZ} 
depend on the choice of the positive integers $M=(M_1,\dots,M_n)$ in congruences \Ref{M_i}. 
Given $M$ satisfying \Ref{M_i}, denote by
\bean
\label{Def M_M}
\mathcal{M}_{M}\,=\,\Big\{ \sum_{l} c_l(z) P^{lp-1}(z,M) \ |\ c_l(z)\in\F_p[z^p]\Big\},
\eean
the $\F_p[z^p]$-module generated by the solutions $P^{lp-1}(z,M)$. 

\vsk.2>
Let $M=(M_1,\dots,M_n)$ and $M\rq{}=(M_1\rq{},\dots,M_n\rq{})$ be two vectors of positive integers, each
satisfying congruences \Ref{M_i}. We say that $M\rq{} > M$ if $M_i\rq{}\geq M_i$ for $i=1,\dots,n$ and there exists $i$ such that
$M_i\rq{}>M_i$.

\vsk.2>
Assume that $M\rq{}>M$. Then  $M_i\rq{} = {M}_i + p N_i$ for some $N_i \in \mathbb{Z}_{\geq0}$.
  We have
\bean
\label{MbM}
P(t,z,M\rq{}) = \Big(\prod_{i=1}^n(t-z_i)^{p N_i}\Big) P(t,z,M) =
\Big(\prod_{i=1}^n(t^p-z_i^p)^{N_i}\Big) P(t,z,M)\,.
\eean
This identity defines an embedding  of modules,
\bean
\label{emb}
\phi_{M\rq{},M} :\mc M_{M\rq{}} \hookrightarrow \mc M_{M}\,.
\eean
Namely, formula \Ref{MbM}, allows us to present any solution $P^{l\rq{}p-1}(z,M\rq{})$, coming
 from  the Taylor expansion of the left-hand side,
 as a linear combination 
of the solutions $P^{lp-1}(z,M)$, coming from the Taylor expansion of
the right-hand side,
 with coefficients in $\F_p[z^p]$.
\vsk.2>

Clearly, if $M\rq{}\rq{}>M\rq{}$ and $M\rq{}>M$, then $M\rq{}\rq{}>M$ and
\bean
\label{comp}
\phi_{M\rq{},M}\,\phi_{M\rq{}\rq{},M\rq{}}\,=\,\phi_{M\rq{}\rq{},M}\,.
\eean

\begin{thm}
\label{thm iso}
 For any $M\rq{}>M$, the embedding
$\phi_{M\rq{},M} :\mc M_{M\rq{}} \hookrightarrow \mc M_{M}$ is an isomorphism.
\end{thm}

\begin{proof}

Let $v\rq{}$ and ${v}$ be the greatest integers such that $v\rq{}p-1 \leq \deg_t P (t,z,M\rq{})$ 
and ${v}p-1 \leq \deg_t  P(t,z,M)$. 
Comparing the coefficients in \Ref{MbM}, we observe that
\begin{gather}
\label{P=matrixQ}
     \begin{bmatrix}
       P^{v\rq{}p-1}(z,M\rq{}) \\
       \vdots \\
       \vdots \\ 
       P^{(v\rq{}-{v}+1)p-1}(z,M\rq{}) \\
     \end{bmatrix} =
                \begin{bmatrix}
                 1 & 0  & \cdots & 0 \\
                 * & \ddots & \ddots & \vdots \\
                 \vdots & \ddots & \ddots & 0 \\
                 * & \cdots & * & 1  \\
                \end{bmatrix} \cdot
                                \begin{bmatrix}
                                P^{{v}p-1}(z, M) \\
                                \vdots \\
                                \vdots \\ 
                                P^{p-1}(z,M) \\
                                \end{bmatrix},
\end{gather}
where all the diagonal entries are $1$'s and stars denote some polynomials in $\FFF_p[z^p]$.
Hence this matrix is invertible and therefore 
$\mc M_M\subset \mc M_{ M\rq{}}$.  The theorem is proved.
\end{proof}

The set of tuples $M=(M_1,\dots,M_n)$ satisfying \Ref{M_i} has the minimal element
$\bar M=(\bar M_1,\dots,\bar M_n)$, where 
$\bar{M}_i$ is  the minimal positive integer satisfying \Ref{M_i}.
Hence for any $M=(M_1,\dots,M_n)$ satisfying \Ref{M_i} we have an
isomorphism 
\bean
\label{MbarM}
\phi_{M,\bar M} :\mc M_{M} \hookrightarrow \mc M_{\bar M}\,.
\eean

\vsk.2>

The module $\mc M_M$, which does not depend on the choice of $M$, will be called
the {\it module of arithmetic solutions} and 
 denoted by $ \mc M_\La$, where $\La$ see in Section \ref{sec ass}.

\subsection{The module  is free}
\label{sec free}

Recall $\La=(\La_1,\dots,\La_n)\in\Z_{>0}^n$ with $\La_i<q$ for $i=1,\dots,n$. Recall
 $\bar M =(\bar M_1,\dots, \bar M_n)$ defined in Section \ref{sec dM}.
Denote
\bean
\label{def dL}
d_\La\,:=\,\Big[ \sum_{i=1}^n\bar M_i / p \Big],
\eean
the integer part of the number $ \sum_{i=1}^n\bar M_i / p$.

\vsk.2>
Consider the module $\mc M_{\bar M}$ spanned over $\F_p[z^p]$
by the solutions $P^{lp-1}(z, \bar M)$, corresponding to $\bar M$.
The range for the index $l$ is defined by  the inequalities 
 $0< lp-1\leq \sum_{i=1}^n \bar M_i -1$. This means that
$l=1,\dots, d_{\La}$.

\begin{thm}
\label{thm lin ind}
The solutions  $P^{lp-1}(z, \bar M)$, $l=1,\dots, d_\La$,
 are linearly independent over the ring $\mathbb{F}_p [z^p]$,
that  is, if \ $\sum_{l=1}^{d_\La} c_l(z) P^{lp-1}(z, \bar M) = 0$ for some 
$c_l (z)  \in \mathbb{F}_p [z^p]$, then $c_l (z) = 0$ for all $l$.

\end{thm}

\begin{cor}
\label{cor 3.3}
The module $\mc M_\La$ of arithmetic solutions is  free  of rank $d_\La$.

\end{cor}

For $q=2$ this theorem is Theorem 3.1 in \cite{V5}. The proof of  Theorem \ref{thm lin ind} below
is the same as the proof of 
 \cite[Theorem 3.1]{V5}.

\begin{proof}

For $l=1,\dots,d_\La$, the coordinates of the vector 
\bea
P^{lp-1}(z,\bar M)
=(P_1^{lp-1}(z,\bar M),\dots,P^{lp-1}_n(z,\bar M))
\eea
 are homogeneous polynomials in $z_1, \dots, z_n$ of
 degree  $\sum_{i=1}^n\bar M_i - lp$ and
\bea
P^{lp-1}_j(z,\bar M)= \sum P^{lp-1}_{j; \ell_1,\dots, \ell_n} z_1^{\ell_1} \dots z_n^{\ell_n},
\eea 
where the sum is over the elements of the set 
\bea
\Ga^{l}_j = \{ (\ell_1, \dots, \ell_n) \in \Z^n_{\geq 0} \, | \, \sum_{s=1}^n \ell_s = \sum_{i=1}^n\bar M_i - lp, 
\, 0 \leq \ell_j \leq \bar M_j -1, \ 0 \leq \ell_i \leq \bar M_i \, \on{for} \, i \ne j \}  
\eea
and
\bea
 P^{lp-1}_{j; \ell_1,\dots, \ell_n} = (-1)^{\sum_{i=1}^n\bar M_i - lp}
\binom{\bar M_j -1}{\ell_j}\prod_{i\ne j}\binom{\bar M_i}{\ell_i}  \ \in \ \FFF_p.
\eea
Notice that all coefficients $ P^{lp-1}_{j; \ell_1,\dots, \ell_n}  $ are nonzero.
Hence each solution $ P^{lp-1}(z,\bar M)$ is nonzero.

\vsk.2>

We show that the first coordinates  $P_1^{lp-1}(z,\bar M)$, $l=1,\dots, d_\La$, 
are linearly independent over the ring $\FFF_p[z^p]$.

\vsk.2>
Let $\bar{\Ga}^{l}_1 \subset \FFF_p^n$ be the image of the
set $\Ga^{{l}}_1$ under the natural projection $\Z^n\to \FFF_p^n$.
The points of $\bar{\Ga}^{{l}}_1$ are in bijective correspondence with the points of
$\Ga^{{l}}_1$, since $\bar M_i < p$ for all $i$. 
Any two sets $\bar{\Ga}^{{l}}_1$ and $\bar{\Ga}^{l\rq{}}_1$ do not intersect, if ${l} \ne {l}'$.

For any ${l}$ and any nonzero polynomial $c_l(z)\in \FFF_p[z^p]$, consider
the nonzero polynomial $c_l(z) P^{{lp-1}}_1(z,\bar M)\in \FFF_p[z_1,\dots,z_n]$ and the set
$\Ga^{{l}}_{1,c_l}$ of vectors $\vec\ell \in\Z^n$ such that the monomial
$z_1^{\ell_1}\dots z_n^{\ell_n}$ enters $c_l(z) P^{lp-1}_1(z,\bar M)$ with nonzero coefficient.
Then the natural projection of $\Ga^{{l}}_{1,c_l}$ to $\FFF_p^n$
coincides with $\bar \Ga^{{l}}_1$. Hence the polynomials $P^{{l}}_1(z,\bar M)$,  $l=1,\dots, d_\La$,
are linearly independent over the ring $\FFF_p[z^p]$.
\end{proof}

\subsection{Fusion  of arithmetic solutions} 
\label{sec fus}

Consider solutions $I(z_1,\dots,z_n)$ of the 
KZ differential equations over $\C$ with values in some tensor product $V_1\otimes\dots\ox V_n$.
Assume that $z_1,\dots, z_n$ tend to some limit $\tilde z_1,\dots,\tilde z_{\tilde n}$, in which some groups of
the points 
$z_1,\dots,z_n$ collide.
It is well-known that under this limit the leading term of asymptotics of solutions  satisfies
 the KZ equations with respect to $\tilde z_1,\dots,\tilde z_{\tilde n}$ with values in the tensor product
$\tilde V_1\otimes \dots\ox\tilde V_{\tilde n}$, where each $\tilde V_j$ is the tensor product of some of
$V_1,\dots,V_n$.

In this section we show that the arithmetic solutions of the KZ equations  have a similar  functorial property but even simpler.

\vsk.2>

Namely, let  $\La=(\La_1,\dots,\La_n)\in\Z_{>0}^n$ with $\La_i<q$ for $i=1,\dots,n$,
and $\tilde \La=(\tilde  \La_1,\dots,\tilde \La_{\tilde n})\in\Z_{>0}^{\tilde n}$
with $\tilde \La_j<q$ for $j=1,\dots,\tilde n$.

\vsk.2>
Assume that there is a partition $\{1,\dots,n\} = I_1\cup\dots \cup I_{\tilde n}$ such that
\bean
\label{fus m}
\sum_{i\in I_j}\La_i = \tilde \La_j\,,\qquad j=1,\dots,\tilde n\,.
\eean
We say that $\tilde \La$ is a {\it fusion} of $\La$.

\vsk.2>
Recall the modules  of arithmetic solutions
$\mc M_\La$ and $\mc M_{\tilde \La}$. We define an epimorphism of modules,
\bean
\label{fu mor}
\psi_{\La,\tilde \La}\,:\, \mc  M_\La\,\to\, \mc M_{\tilde \La}\,,
\eean
as follows.

Let $M=(M_1,\dots, M_n)$ be a vector of positive integers such that
$M_i\equiv -\La_i/q $ mod $p$ for $i=1,\dots,n$.
Define  $\tilde M=(\tilde M_1,\dots,\tilde M_{\tilde n})$ by the formula
\bea
\tilde M_j = \sum_{i\in I_j} M_i\,,\qquad j=1,\dots,\tilde n\,.
\eea
Then 
$\tilde M_j\equiv -\tilde \La_j/q $ (mod $p$) for $j=1,\dots,\tilde n$.

Consider the module $\mc M_{M}$ generated by the solutions
$(P^{lp-1}(z,M))_l$, which are $n$-vectors of polynomials in $z_1,\dots,z_n$.
Consider the module $\mc M_{\tilde M}$ generated by the solutions
$(P^{lp-1}(\tilde z, \tilde M))_l$, which are $\tilde n$-vectors of polynomials 
in $\tilde z_1,\dots,\tilde z_{\tilde n}$.

Define the module homomorphism 
\bea
 \psi_{\La,\tilde \La}:\mc M_\La = \mc  M_{M}\,\to\, \mc M_{\tilde M}=\mc M_{\tilde \La}\,,
\eea
as the map which sends a generator $P^{lp-1}(z,M)$ to the generator $P^{lp-1}(\tilde z,\tilde M)$.

\vsk.2>
On the level of coordinates of these vectors,  the vectors $P^{lp-1}(z,M)$ and $P^{lp-1}(\tilde z,\tilde M)$
have the following relation. Choose a coordinate
$P^{lp-1}_a(z,M)$ of $P^{lp-1}(z,M)$,  replace in it every $z_i$ with $\tilde z_j$ if $i\in I_j$,
then the resulting polynomial $P^{lp-1}_a(z(\tilde z),M)$
 equals the $b$-th coordinate $P^{lp-1}_b(\tilde z,\tilde M)$ of
$P^{lp-1}(\tilde z,\tilde M)$ if $a\in I_b$.

It is easy to see that the homomorphism $ \mc  M_{\La}\,\to\, \mc M_{\tilde \La}\,,$ does not depend of the choice
of $M$ solving congruences \Ref{M_i}.

\vsk.2>
Also it is easy to see that if $\tilde \La$ is a fusion of $\La$ and $\hat \La$ is a fusion of $\tilde \La$, then
$\hat \La$ is a fusion of $\La$ and 
\bean
\label{comp f}
\psi_{\tilde \La, \hat \La}\,\psi_{\La, \tilde \La}\,=\,\psi_{\La,\hat \La}\,.
\eean

\section{Binomial coefficients modulo $p$}
\label{sec bin coef}

\subsection{Lucas' theorem}

\begin{thm} [\cite{L}] 
\label{thm L}

For nonnegative integers $m$ and $n$ and a prime $p$,
 the following congruence relation holds:
\bean
\label{Lucas}
\binom{n}{m}\equiv \prod _{i=0}^k \binom{n_i}{m_i}\quad (\on{mod}\ p),
\eean
where $m=m_{k}p^{k}+m_{k-1}p^{k-1}+\cdots +m_{1}p+m_{0}$ 
and  $n=n_{k}p^{k}+n_{k-1}p^{k-1}+\cdots +n_{1}p+n_{0}$
are the base $p$ expansions of $m$ and $n$ respectively. This uses the convention that 
$\binom{n}{m}=0$ if $n<m$.
\qed
\end{thm}

On Lucas\rq{} theorem see for example, \cite{Me}.

\subsection{Factorization of  $\binom{-1/q}{m}$ modulo $p$}
\label{ssec fact} 

In the next sections we will use the binomial coefficients $\binom{-1/q}{m}$ modulo $p$.
Recall that $p>q$ are prime numbers. Denote by $d$ the order of $p$ modulo $q$, that is the least integer $k$ such that
$p^k\equiv 1$ (mod $q$).

Let 
\bean \label{p exp 1}
\frac{p^d-1}{q} = A_0 + A_1 p + A_2 p^2 + \dots + A_{d-1} p^{d-1}
\eean
be the base $p$ expansion of $(p^d-1)/q$.

\begin{thm} 
\label{thm 1/q}

Let $m$ be a positive integer with the base $p$ expansion given by
\bean
 \label{p exp m}
m=m_b p^b + m_{b-1}p^{b-1}+ \dots + m_1 p + m_0.
\eean 
Then the binomial coefficient $\binom{-1/q}{m}$ is well-defined modulo $p$ and the following congruence holds
\bean
\label{-1/q m}
\binom{-1/q}{m} \equiv  \prod_{j\geq 0} 
\prod_{i=0}^{d-1} \binom{A_i}{m_{dj+i}} \quad (\on{mod}\ p).
\eean 
\end{thm}

The case $q=2$ was considered in \cite{V5}.

\begin{proof}
The proof of this theorem is based on Lucas\rq{} theorem and the following three lemmas.

\begin{lem}
\label{lem aux}
For any positive integer $c$ 
\bean
\label{p exp 2}
\frac{p^{dc}-1}{q} = 
\sum_{j=0}^{c-1} \Big{(} \sum_{i=0}^{d-1} A_i p^i \Big{)} p^{jd}.
\eean
\end{lem}
\begin{proof}
We have 
\bea
\frac{p^{dc}-1}{q} = \frac{(p^d-1)(p^{d(c-1)}+\dots+p+1)}{q}
=\sum_{j=0}^{c-1} \Big{(} \sum_{i=0}^{d-1} A_i p^i \Big{)} p^{jd},
\eea
where we use \Ref{p exp 1} to obtain the last equality. 
\end{proof}

\begin{lem} \label{lem fact}
Let $m$ be a positive integer with the base $p$ expansion given by \Ref{p exp m}.
If $c$ is a positive integer such that $dc>b$, then 
\bean \label{fact}
\binom{(p^{dc}-1)/q}{m} \equiv \prod_{j=0}^{c-1} 
\prod_{i=0}^{d-1} \binom{A_i}{m_{dj+i}} \quad (\on{mod}\ p).
\eean 
\end{lem}
\begin{proof}
The lemma follows from Theorem \ref{thm L} and Lemma \ref{lem aux}.
\end{proof}

Given a prime $p$, define the
{\it  $p$-adic norm} on $\Q$ as follows. 
Any nonzero rational number $x$ can be represented uniquely 
by $x=p^{\ell}(r/s)$, where $r$ and $s$ are integers not 
divisible by $p$. Set $|x|_p=p^{-{\ell}}$. Also define the 
$p$-adic value $|0|_p=0$. We call $p^{\ell} (r/s)$ the \emph{$p$-reduced presentation} of $x$.

\begin{lem} \label{lem wd}
Let $m$ be a positive integer with the base $p$ expansion
 given by \Ref{p exp m}. Then the following statements hold.
\begin{enumerate}
\item[(i)] The binomial coefficient $\binom{-1/q}{m}$ is well-defined modulo $p$.
\item[(ii)] For any positive integer $c$ such that $dc > b+1$ we have 
\bean
\label{equiv}
\binom{-1/q}{m} \equiv \binom{(p^{dc}-1)/q}{m} \qquad (\on{mod}\ p).
\eean
\end{enumerate}
\end{lem}

\begin{proof}
We have
\bean
\label{pr1}
\binom{-1/q}{m} = (-1/q)^m \prod_{\ell=0}^{m-1} \frac{q\ell+1}{\ell+1}.
\eean
We also have 
\bean
\label{pr2}
\binom{(p^{dc}-1)/q}{m} = (-1/q)^m \prod_{\ell=0}^{m-1} \frac{q \ell+1-p^{dc}}{\ell+1} = (-1/q)^m \prod_{\ell=0}^{m-1} \frac{q\ell+1}{\ell+1} + \dots
\eean

For each $\ell = 0, \dots, m-1$ we have $\big| \frac{q\ell+1 - p^{dc}}{\ell+1} \big|_p = \big| \frac{q\ell+1}{\ell+1} \big|_p $, since $q\ell+1 < qm < qp^{b+1} <p^{b+2} \leq p^{dc}$. 

By the same reasoning for every $\ell$ the power of $p$ in the $p$-reduced
 presentation of the number $\frac{p^{dc}}{\ell+1}$ is greater than the 
power of $p$ in the $p$-reduced presentation of $\frac{q\ell+1}{\ell+1}$. 
This observation implies that in the right-hand side of \Ref{pr2} the power
 of $p$ in the $p$-reduced presentation of terms denoted by '$\dots$'
 is greater than the corresponding power in the $p$-reduced presentation 
of $(-1/q)^m \prod_{\ell=0}^{m-1} \frac{q\ell+1}{\ell+1}$. Since the
 left-hand side of \Ref{pr2} is an integer, we conclude that in the 
$p$-reduced presentation of the binomial coefficient $\binom{-1/q}{m}$ the 
power of $p$ is non-negative, i.e. $\binom{-1/q}{m}$ is well-defined modulo $p$. 
This gives part (i). Moreover, the following congruence holds
\bean
\label{pr3}
\binom{(p^{dc}-1)/q}{m} \equiv (-1/q)^m \prod_{\ell=0}^{m-1} \frac{q\ell+1}{\ell+1} \qquad (\on{mod}\ p).
\eean 
Formulas \Ref{pr1} and \Ref{pr3} give \Ref{equiv}.
\end{proof}
Theorem \ref{thm 1/q} follows from Lemmas \ref{lem fact} and \ref{lem wd}.
\end{proof}

\subsection{Factorization of $\binom{u/v}{m}$ modulo $p$}

\label{ssec g.fact}

Here are some generalizations of Theorem \ref{thm 1/q}.
Let $p$ be  prime. Let 
$u, v$  be relatively prime integers with $0 < u ,v < p$.

\vsk.2>
Assume that there exists a positive integer $\ell$ such that 
\bean
\label{assump}
p^{\ell} + u \,\equiv\, 0 \qquad (\on{mod} \ v)\,.
\eean
Let 
\bean
\label{p exp 4}
\frac{p^\ell+u}{v} = B_0 + B_1 p + B_2 p^2 + \dots + B_{\ell-1} p^{\ell-1}
\eean
be the base $p$ expansion of $(p^\ell+u)/v$\,. 

\vsk.2>
Let $\phi(v)$ be the number of positive divisors of $v$ and 
\bean 
\label{p exp 5}
\frac{p^{\phi(v)}-1}{v} = C_0 + C_1 p + C_2 p^2 + \dots + C_{\phi(v)-1} p^{\phi(v)-1}
\eean
the base $p$ expansion of $(p^{\phi(v)}-1)/v$\,.

Let $m$ be  a positive integer with the base $p$
expansion
\bea
m=m_b p^b + m_{b-1}p^{b-1}+ \dots + m_1 p + m_0\,.
\eea

\begin{thm}
\label{thm u/v}

Under these assumptions
the binomial coefficient $\binom{u/v}{m}$ is well-defined modulo $p$
 and the following congruence holds:
\bean 
\label{u/v m}
\binom{u/v}{m} \equiv \left[ \prod_{i=0}^{\ell-1} \binom{B_i}{m_i} \right] 
\cdot \left[ \prod_{j \geq 0}
 \prod_{i=0}^{\phi(v)-1} \binom{C_i}{m_{\ell+j \phi(v) + i}} \right] \qquad (\on{mod} \ p)\,.
\eean
\end{thm}

\begin{proof}
The proof of Theorem \ref{thm u/v} is parallel to the proof of Theorem \ref{thm 1/q}
and follows from the analysis of the base $p$ expansion of the integers of the form
\bea
\frac{p^{k\phi(v) +  \ell}+ u}v\, =\, p^\ell\,\frac{p^{k\phi(v)}-1}v\,+\,\frac{p^\ell+ u}v
\eea
 for $k\in\Z_{\geq 0}$.
\end{proof}

In the same way we  prove the following statement.

\begin{thm}
\label{thm gen}
Let $p$ be  a  prime. Let $u, v$  be integers such that $p, u, v$
are pairwise  relatively prime.
Let
\bean
\label{gen u/v}
\frac uv\,=\, A_0 + A_1p + A_2p^2+\dots\,,
\qquad  0\leq A_i< p \quad\on{for}\ i\in\Z_{\geq 0}\,,
\eean
be the $p$-adic presentation of $u/v$.
Let $m$ be  a positive integer with the base $p$
expansion
\bea
m=m_b p^b + m_{b-1}p^{b-1}+ \dots + m_1 p + m_0\,.
\eea
Then
\bean
\label{Lu gen}
\binom{u/v}{m} \equiv \prod_{i\geq 0} \binom{A_i}{m_i} 
 \qquad (\on{mod} \ p)\,.
\eean
\qed
 
\end{thm}

\subsection{Map $\eta$.}\label{ssec eta}

Recall that  $p$ and $q$ be prime numbers, $p>q$, and  $d$ the order of $p$ modulo $q$.
Define the map
\bean
\label{eta}
\eta\ :\ \{1,\dots,q-1\} \ \to\ \{1,\dots,q-1\}\,, \qquad a \ \mapsto\ \eta(a)\,, 
\eean
the division by $p$ modulo $q$. More precisely,  $\eta(a)$ is defined by the conditions
\bean
\label{congr}
\eta(a) p \,\equiv\, a \quad (\on{mod}\ q), \qquad 1\leq \eta(a) < q\,.
\eean

\vsk.2>
Denote  $\eta^{(s)}=\eta \circ \eta \circ \dots \circ \eta$ the $s$-th iteration of $\eta$.
We have  $\eta^{(s +d)} = \eta^{(s)}$.
Define
\bean
\label{a_s}
a_s \,:= \,\eta^{(s)}(1)\,.
\eean
We have $a_{s+d}=a_s$ and $a_0=1$.
\vsk.2>

The integer $a_1$ will play a special role in the next sections. It is the unique integer such that 
\bean
\label{a_1}
1\leq a_1< q\quad \on{and}\quad  q\,\vert \,(a_1p-1)\,.
\eean

\vsk.2>

\begin{exmp}
\label{ex a_s}
Let $p=5, \, q=3$. The order $d$  of 5 modulo 3 is $2$. The map $\eta : \{1,2\} 
\to \{2,1\}$ is the transposition. We have  $a_0 = 1$, $a_1=2$, $a_2=1$.
\end{exmp}

\subsection{Formulas for $A_s$}\label{ssec A_s}
Recall 
the base $p$ expansion
\bean
\label{p exp 3} 
\frac{p^d-1}{q} = A_0 + A_1 p + A_2 p^2 + \dots + A_{d-1} p^{d-1}\,.
\eean

\begin{lem}\label{A and a}
The integers $A_s$ are given by the formula 
\bean
\label{A_s}
A_s = \frac{a_{s+1}p-a_s}{q}\,, \quad s =0,\dots,d-1\,.
\eean
Moreover, $A_s > 0$ for $s=0,\dots,d-1$.
\end{lem}
\begin{proof}
We have
\bean
\label{A_s est} 
0 < \frac{a_{s+1}p-a_s}{q}  < p 
\eean
for every $s$. Indeed, since $1 \leq a_s \leq q-1$ and $p>q$, we have
\bea
0 < \frac{p-a_s}{q} \leq \frac{a_{s+1}p-a_s}{q} < \frac{a_{s+1}p}{q} < p\,.
\eea
We show that \Ref{p exp 3} holds for $A_s$ given by \Ref{A_s}. Indeed
\bea
q ( A_0 + A_1 p + \dots + A_{d-1} p^{d-1} ) 
= q \left( \frac{a_1p-1}{q} + \frac{a_2p-a_1}{q} p + \dots + \frac{p-a_{d-1}}{q} p^{d-1} \right)
= p^d-1\,.
\eea
\end{proof}

\section{Solutions  for the case $\La_i=1$, $i=1,\dots,n$}
\label{sec m=1}

In Sections \ref{sec m=1}, \ref{sec Cartier-Manin},  and \ref{sec 8} we  study 
the arithmetic solutions of the KZ equations over $\F_p$ in the case $\La_i=1$, $i=1,\dots,n$.

\subsection{Basis of arithmetic solutions}
\label{sec Bas}

Recall that $p, q$ are prime numbers, $p> q$.
 Assume that
\bean
\label{p q k}
n\,=\,kq \,+\,1\
\eean 
for some positive integer $k$ and
\bean
\label{all m_i=1}
\La = (1,\dots,1) \in \Z^n\,.
\eean

\vsk.2>
The minimal positive solution of the system of the congruences
\bean
\label{M_i m_i=1}
M_i\equiv -1/q
\quad 
(\on{mod} \,p)\,, \qquad i=1,\dots,n\,,
\eean
is the vector 
\bean
\label{baR M}
\bar M =(\bar M_1,\dots, \bar M_n) := ((a_1p-1)/q, \dots, (a_1p-1)/q)\,,
\eean
where $a_1$ is introduced in \Ref{a_1}.

\vsk.2>
Recall that the rank of the module of arithmetic solutions $\mc M_\La$ in this case equals
\bean
\label{d m1}
d_\La\, =\,\left[ n (a_1p-1)/qp   \right]\,,
\eean
see Corollary \ref{cor 3.3}.

\begin{lem}
\label{lem d m1}
If $p>n$, then $d_\La = a_1 k$.

\end{lem}

\begin{proof} We have
\bea 
d_\La =\left[ n \frac{a_1p-1}{qp}\right] = \left[ \frac{na_1}q-\frac n{qp}\right]
=
 \left[ ka_1 +\frac{a_1}q - \frac n{pq}\right].
\eea
Notice that   $1/q\leq a_1/q <1$ and $n/qp< 1/q$. Hence $d_\La=a_1k$.
\end{proof}

\vsk.2>
The master polynomial is 
\bean
\label{mast p}
\Phi_p(x,z,\bar M) = \prod_{i=1}^n(x-z_i)^{(a_1p-1)/q} \,.
\eean
The $n$-vector of polynomials $P(x,z,\bar M)$ is
\bean
\label{P(x,z,M)}
\phantom{aaa}
P(x,z, \bar M)=(P_1(x,z, \bar M),\dots,P_n(x,z, \bar M)),\quad P_j(x,z, \bar M)=  \frac {\Phi_p(x,z, \bar M)}{x-z_j}\,,
\eean
with the Taylor expansion
\bean
\label{te}
\phantom{aaaaaa}
P(x,z, \bar M)= \sum_{i=0}^{n(a_1p-1)/q -1}  P^{i}(z, \bar M)\, x^i, 
\quad P^{i}(z, \bar M)=(P^{i}_1(z,\bar M),\dots,P^{i}_n(z,\bar M))\,.
\eean
The coefficients $P^{i}(z, \bar M)$  with $i=lp-1$ are  solutions of \Ref{KZ}.
There are  $ka_1$ of them. They are linearly independent by Theorem \ref{thm lin ind}.
Denote them 
\bean
\label{Imz}
I^m (z) = ( I^m_1 (z), \dots, I^m_n (z) )\,, \qquad
I^m (z) := P^{(a_1k-m)p+p-1} (z,\bar M)\,,
\eean
for $m=1,\dots, a_1k$.
The coordinates of the $n$-vector of polynomials $I^m (z)$ are {\it homogeneous polynomials in $z$
of degree }
\bea
(a_1p-1)/q+(m-1)p-k\,.
\eea
For example, $I^1=P^{a_1kp-1} (z,\bar M)$ is a homogeneous polynomial in $z$ of degree 
$(a_1p-1)/q-k$, and $I^{a_1k}=P^{p-1} (z,\bar M)$
is a homogeneous polynomial in $z$ of degree 
$(a_1p-1)/q-k+(a_1k-1)p$.

\subsection{Change of the basis of $\mc M_\La$}
For future use we introduce a new basis of $\mc M_\La$.
For $m=1,\dots,a_1k$, define
\bean
\label{J&I}
J^m(z)=\sum_{l=1}^m  I^{m+1-l}(z) \, z_1^{(l-1)p} \,  \binom{a_1k-m -1+ l}{a_1k-m},
\eean
that is, 
\bea
J^1(z) &=& I^1(z),
\\
J^2(z) &=&  I^{1}(z) \, z_1^{p} \, \binom{a_1k-1}{a_1k-2}  + I^2(z),
\\
J^3(z) &=&  I^{1}(z) \, z_1^{2p} \, \binom{a_1k-1}{a_1k-3}  +  I^{2}(z) \, z_1^{p} \, \binom{a_1k-2}{a_1k-3}  + I^3(z), 
\eea
and so on.
The coordinates of the $n$-vector of polynomials $J^m (z)$ are {\it homogeneous polynomials in $z$
of degree }
\bea
(a_1p-1)/p+(m-1)p-k\,.
\eea

\begin{lem}
The $n$-vectors of polynomials $J^m(z)$ belong to the module $\mc M_\La$ of arithmetic solutions and
form a basis of $\mc M_\La$.
\qed
\end{lem}

The solutions $J^m(z)$ of the KZ equations \Ref{KZ} can be defined as 
coefficients of the Taylor expansion of a 
suitable $n$-vector of polynomials 
similarly to the definition of solutions $I^m(z)$ in \Ref{Imz}.

 Namely,  change  $x$ to $x+z_1$ in  \Ref{mast p} to obtain the polynomial 
\bean
\label{mast p til}
\tilde \Phi_p(x,z) = x^{(a_1p-1)/q} \prod_{i=2}^n (x+z_1-z_i)^{(a_1p-1)/q}  \ \in \ \FFF_p[x,z]\,.
\eean
Consider the associated $n$-vector of polynomials 
\bean
\label{til P}
\tilde P (x,z) = \tilde \Phi_p(x,z) \left( \frac 1x, \frac 1{x+z_1-z_2},\dots, \frac 1{x+z_1-z_n} \right)
\eean
and its Taylor expansion
\bean
\label{te til}
\tilde P(x,z) = \sum_{i=0}^{{n(a_1p-1)/q -1}}  \tilde P^i(z) x^i\,,\qquad \tilde P^i(z) \in \FFF_p[z]^n\,.
\eean

\begin{lem} [{\cite[Lemma 5.2]{V5}}]
\label{lem tP=J}
For $m=1,\dots,a_1k$, we have
\bean
\label{Jmz}
J^m(z) = \tilde P^{(a_1k-m)p+p-1}(z)\,.
\eean
\qed
\end{lem}

\subsection{Change of variables in $J^m(z)$}
\begin{lem}
\label{cor z_j-z_1}
For $m=1,\dots,a_1k$,  the homogeneous polynomial  $J^m(z)$ 
can be written  as a homogeneous
polynomial in  variables $z_2-z_1$, $z_3-z_1$, \dots, $z_n-z_1$, see
\Ref{mast p til} and \Ref{til P}.
\qed
\end{lem}

Hence we may pull out from the polynomial $J^m(z)$ the factor $(z_2-z_1)^{(a_1p-1)/q + (m-1)p-k}$ 
and present $J^m(z)$ in the form
\bean
\label{J&K}
J^m(z) = (z_2-z_1)^{(a_1p-1)/q + (m-1)p-k}\, K^m(\la)\,,
\eean
where $\la: =  (\la_3,\dots,\la_n)$,
\bean
\label{la_j}
\la_i\,:=\,\frac{z_i-z_1}{z_2-z_1}, \qquad i = 3,\dots,n,
\eean
and $K^m(\la)$ is a suitable  $n$-vector of polynomials in $\la$.

\vsk.2>
Another way to define  the $n$-vector $K^m(\la)$ is as follows.

\vsk.2> 

Consider the polynomial 
\bean
\label{mast p hat}
\hat \Phi_p (x,\la) = x^{(a_1p-1)/q}(x-1)^{(a_1p-1)/q} \prod_{i=3}^n (x-\la_i)^{(a_1p-1)/q} \ \in  \ \FFF_p[x,\la]
\eean
and the associated $n$-vector of polynomials 
\bean
\label{hat P}
\hat P(x,\la) = \hat \Phi_p (x,\la) \left( \frac 1x, \frac 1{x-1}, \frac 1{x-\la_3}, \dots, \frac 1{x-\la_n} \right)
\eean
with Taylor expansion 
\bean
\label{te hat}
\hat P(x,\la) = {\sum}_{i=0}^{n(a_1p-1)/q -1}  \hat P^i(\la) x^i\,, \qquad
\hat P^i(\la) \in \FFF_p[\la]^n\,.
\eean

\begin{lem}
 For $m=1,\dots, a_1k$, we have
\bean
\label{Kmla}
K^m(\la) = \hat P^{(a_1k-m)p+p-1} (\la)\,,
\eean
cf. formulas \Ref{Imz} and \Ref{Jmz}.
\qed
\end{lem}

\subsection{Formula for $K^m(\la)$}\label{ssec fmla K}

For $m=1, \dots, a_1k,$ denote
\bean
\label{Del fh}
\qquad  \Delta^{m } = \Big\{(\ell_3,\dots,\ell_n)\in \Z^{n-2}_{\geq 0} \ | \  
0 \leq \sum_{i=3}^n \ell_i + k -(m-1)p \leq (a_1p-1)/q, 
\\
\phantom{aaaaaaaaaaaaaaaaaa} \ell_j \leq (a_1p-1)/q \ \on{for} \ j=3,\dots,n \Big\}. \phantom{aaaaa} \nonumber
\eean

\begin{thm}\label{thm K}
For $m=1,\dots,a_1k$, we have
\bean
\label{K=}
K^m(\la)=\sum_{(\ell_3,\dots,\ell_n) \in\Delta^m} K^m_{\ell_3,\dots,\ell_n} (\la)\, ,
\eean
where 
\bean
\label{K ell}
&&
 K^m_{\ell_3,\dots,\ell_n} 
(\la) = (-1)^{(a_1p-1)/q +(m-1)p-k} \binom{(a_1p-1)/q}{\sum_{i=3}^{n}\ell_i +k-(m-1)p} 
\\
\notag
&&
\phantom{a}
\times\,
\prod_{i=3}^{n} \binom{(a_1p-1)/q}{\ell_{i}} \ \la_3^{\ell_3}\dots\la_n^{\ell_n} \,
\Big(1, 
-q \Big( \sum_{i=3}^n \ell_i + k \Big), q\ell_3+1, \dots, q\ell_n+1\Big) 
\eean
modulo $p$. In particular, all coefficients  $ K^m_{\ell_3,\dots,\ell_n}(\la) $ are nonzero. 
\end{thm}

\begin{proof} 
We have
\bea
K_1^m(\la) = (-1)^{(a_1p-1)/q + (m-1)p -k}{\sum_{\Delta}} \,\binom{(a_1p-1)/q}{\ell_2}
\dots \binom{(a_1p-1)/q}{\ell_n} \la_3^{\ell_3}\dots\la_n^{\ell_n},
\eea
where 
\bea
\Delta &=& \{(\ell_2,\dots,\ell_n)\in \Z^{n-1}_{\geq 0} \ | \  \sum_{i=2}^n \ell_i =(a_1p-1)/q + (m-1)p-k,\ \ \ell_i \leq (a_1p-1)/ q\}.
\eea
Expressing $\ell_2$ from the conditions defining $ \Delta $ we  write 
\bea
&&
K^m_1(\la) 
=
 (-1)^{(a_1p-1)/q + (m-1)p -k} 
\\
&&
\phantom{aaa}
\times\,
{\sum_{\Delta^m}} \,\binom{(a_1p-1)/q}{{\sum}_{i=3}^{n}\ell_i - (m-1)p+k} \prod_{i=3}^n \binom{(a_1p-1)/q}{\ell_i} \la_3^{\ell_3}\dots\la_n^{\ell_n}.
\eea
Similarly we have
\bea
K_2^m(\la) \ = \ (-1)^{(a_1p-1)/q+(m-1)p-k} 
 {\sum_{\Delta'}} \binom{(a_1p-1)/q-1}{\ell_2} \prod_{i=3}^n \binom{(a_1p-1)/q}{\ell_i} \la_3^{\ell_3}\dots\la_n^{\ell_n},
\eea
where
\bea
\Delta' &=& \{(\ell_2,\dots,\ell_n)\in \Z^{n-1}_{\geq 0} \ | \  \sum_{i=2}^n \ell_i =  (a_1p-1)/q + (m-1)p-k,\ \ 
\\
&&
\phantom{aaaaaa}
\ell_2 \leq (a_1p-1)/q -1 \ \on{and}\ 
\ell_i \leq (a_1p-1)/q \ \on{for}\  i > 2\}.
\eea
Expressing $\ell_2$ from the conditions defining $\Delta'$ we write 
\bea
K_2^m(\la) &= & (-1)^{(a_1p-1)/q+(m-1)p-k} 
\\
&&
\phantom{a}
\times\ {\sum_{\Delta^m}} \binom{(a_1p-1)/q-1}{\sum_{i=3}^n \ell_i - (m-1)p+k-1} \prod_{i=3}^n \binom{(a_1p-1)/q}{\ell_i} \la_3^{\ell_3}\dots\la_n^{\ell_n} \phantom{aaaaa}
\\
&=&
 (-1)^{(a_1p-1)/q+(m-1)p-k} \ \frac{\sum_{i=3}^n \ell_i - (m-1)p+k}{(a_1p-1)/q}  \phantom{aaaaaaaaaaaaaa} 
\\
&&
\phantom{a}
\times \ {\sum_{\Delta^m}} \binom{(a_1p-1)/q}{\sum_{i=3}^n \ell_i - (m-1)p+k} \prod_{i=3}^n \binom{(a_1p-1)/q}{\ell_i} \la_3^{\ell_3}\dots\la_n^{\ell_n}\,. 
\eea
For $j=3,\dots,n$ we have
\bea
K_j^m(\la) = (-1)^{(a_1p-1)/q+(m-1)p-k}
 {\sum_{\Delta''}} \binom{(a_1p-1)/q-1}{\ell_j} 
\prod_{i=2 \atop i\ne j}^n \binom{(a_1p-1)/q}{\ell_{i}} \la_3^{\ell_3}\dots\la_n^{\ell_n}, 
\eea
where 
\bea
\Delta'' &=& \{(\ell_2,\dots,\ell_n) \in \Z^{n-1}_{\geq 0} \ | \  \sum_{i=2}^n \ell_i = (a_1p-1)/q+ (m-1)p-k,\ \ 
\\
&&
\phantom{aaaa}
\ell_j \leq (a_1p-1)/q-1 \ \on{and}\ \ell_i \leq (a_1p-1)/q \ \on{for}\ i\ne j\}.
\eea
Expressing $\ell_2$ from the conditions defining $\Delta''$ we write 
\bea
&&
K_j^m(\la)
\,=\, (-1)^{(a_1p-1)/q+(m-1)p-k} 
\\
&&
\phantom{aa}
\times \,{\sum_{\Delta^m}} \binom{(a_1p-1)/q}{\sum_{i=3}^n \ell_i - (m-1)p+k} \binom{(a_1p-1)/q-1}{\ell_j} \prod_{i=3 \atop i\ne j}^n \binom{(a_1p-1)/q}{\ell_{i}}
\la_3^{\ell_3}\dots\la_n^{\ell_n} 
\\
&&
\phantom{aa}
=\, (-1)^{(a_1p-1)/q+(m-1)p-k} \ \big(1 - \frac{\ell_j q}{a_1p-1} \big) 
\\
&&
\phantom{aa}
\times \,
{\sum_{\Delta^m}}\, \binom{(a_1p-1)/q}{\sum_{i=3}^n \ell_i - (m-1)p+k} \prod_{i=3}^n \binom{(a_1p-1)/q}{\ell_{i}} \la_3^{\ell_3}\dots\la_n^{\ell_n}\,.
\eea
The congruences
\bea 
\frac{\sum_{i=3}^n \ell_i - (m-1)p+k}{(a_1p-1)/q} 
&\equiv&
 -q\big{(}{\sum}_{i=3}^n \ell_i + k\big{)} \qquad (\on{mod}\ p), 
\\ 
 1 - \frac{\ell_j q}{a_1p-1} 
&\equiv&
 q\ell_j+1 \qquad (\on{mod}\ p)  
\eea
allow us to rewrite $K^m_j(\la), \ j = 2,\dots,n,$ in the form indicated in the theorem.
\end{proof}

\subsection{Example}
\label{K^m}

For $p=5$, $q=3$, $n=4$, the module $\mc M_\La$ of arithmetic solutions is two-dimensional and is generated
by solutions $J^1$ and $J^2$ which are $4$-vectors, whose coefficients are homogeneous polynomials in
$z_1,z_2,z_3,z_4$ of 
degrees 2 and 7, respectively. 

The corresponding $4$-vectors $K^1$ and $K^2$ are given by the formulas
\bea
K^1(\la_3,\la_4) 
&=&
 (3,1,3,3) + (4,1,1,4) \la_3 + (4,1,4,1) \la_4  
\\
&+&
 (3,3,1,3)\la_3^2 + (4,4,1,1)\la_3 \la_4 + (3,3,3,1) \la_4^2\,,
\\
K^2(\la_3,\la_4) 
&=&
 (2,0,0,3) \la_3^2 \la_4 + (1,0,2,2) \la_3^2 \la_4^2 + (2,0,3,0) \la_3 \la_4^3 
\\
&+&
 (1,2,0,2) \la_3^3 \la_4^2 + (1,2,2,0) \la_3^2 \la_4^3 + (2,3,0,0) \la_3^3 \la_4^3\,.
\eea

\section{Cartier map}
\label{sec Cartier-Manin}

\subsection{Matrices and semilinear algebra, \cite{AH}}
\label{semilin}

\subsubsection{Bases, matrices, and linear operators}

Let $W$ be a vector space over a field $\Bbb K$
with basis $\mathcal{C} = \{w_1,\dots,w_n\}$. Any $w \in W$ is expressible
 as $w = \sum c_i w_i$. Let $[w]_{\mathcal{C}}$ denote the column vector
$[w]_{\mathcal{C}} = (c_1,\dots,c_n)^T$.

Let V be a vector space with basis ${\mathcal{B}} = \{ v_1,\dots, v_m \}$, and  $f: W \to V$  a linear map.
 The matrix of $f$ relative to the 
 bases ${\mathcal{C}}$ and ${\mathcal{B}}$ is 
$[f]_{{\mathcal{B}} \leftarrow {\mathcal{C}}} = (a_{ij})  \in  \on{Mat}_{m\times n} (\mathbb{K})$,
where $f(w_j) = \sum_{i=1}^m a_{ij} v_i$. Matrix multiplication gives
\bea 
[f(w)]_{\mathcal{B}} = [f]_{{\mathcal{B}} \leftarrow {\mathcal{C}}} \ [w]_{\mathcal{C}}.
\eea
 Let $f: V \to V$ be an endomorphism, and  ${\mathcal{B}}$ and ${\mathcal{D}}$ two 
 bases for $V$. Then
\bea
[f]_{\mathcal{D}\leftarrow\mc D} = [\id]_{{\mathcal{D}} \leftarrow {\mathcal{B}}} \ [f]_{\mathcal{B}\leftarrow \mc B}
 \ [\id]_{{\mathcal{B}} \leftarrow {\mathcal{D}}}\,.
\eea
 If $S = [\id]_{{\mathcal{B}} \leftarrow {\mathcal{D}}}$, then 
\bea
[f]_{\mathcal{D}\leftarrow \mc D} = S^{-1} \ [f]_{\mathcal{B}\leftarrow \mc B} \ S\,.
\eea

\subsubsection{Semilinear algebra}

Let $\tau$ be an automorphism of $\mathbb{K}$ and $\si=\tau^{-1}$. 
Then $f:V \to V$ is called {\it $\tau $-linear}, if for $a \in \mathbb{K}$ and $v \in V$, 
\bea 
f(av) = a^{\tau} f(v),
\eea 
where $a^{\tau} = \tau (a)$. 
Let $f(v_j) = \sum_i a_{ij} v_i$. 
If $v= \sum_j c_j v_j$, then
\bea 
f(v) = \sum_j f(c_j v_j) = \sum_j c^{\tau}_j f(v_j) = \sum_j \Big( \sum_i a_{ij} v_i \Big) c^{\tau}_j 
\eea
and so
\bea
[f(v)]_{\mathcal{B}} \,=\, [f]_{\mathcal{B}\leftarrow\mc B} \ [v]^{\tau}_{\mathcal{B}}\,,
\eea
where $B^{\tau}$ is the matrix  obtained by applying $\tau$ to each entry of $B$.

Change of basis is accomplished with $\tau$-twisted conjugacy:
\bea
[f]_{\mathcal{D}\leftarrow \mc D} = [\id]_{{\mathcal{D}} \leftarrow {\mathcal{B}}}
 \ [f]_{\mathcal{B}\leftarrow\mc B} \ [\id]_{{\mathcal{B}} \leftarrow {\mathcal{D}}}^{\tau} 
= S^{-1} \ [f]_{\mathcal{B}\leftarrow\mc B} \ S^{\tau}.
\eea

The iterates of $f$ are represented by 
\bean
\label{iterr} 
[f^{\circ r}] = [f] \ [f]^{\tau} \ [f]^{\tau^2} \dots [f]^{\tau^{r-1}}.
\eean

\subsubsection{Adjoint  map} 
Let $V^*$ be the dual vector space of $V$ and 
 $(\cdot,\cdot)_V : V \times V^* \to \mathbb{K}$ the natural pairing,
 linear with respect to each argument.  Similarly, let
 $W^*$ be the dual vector space  of $W$ and 
 $(\cdot,\cdot)_W: W \times W^* \to \mathbb{K}$ the natural pairing.
 
\vsk.2>

 Let $f: W \to V$ be $\tau$-linear. Define the
{\it adjoint map} $f^*: V^*\to W^*$ by the formula
\bea
(w, f^*(\phi))_W \,=\, (f(w), \phi)_V^\si, \qquad w\in W\,, \ \phi\in V^*\,.
\eea
The map $f^*$ is $\si$-linear,  $f^*(a\phi) = a^\si f^*(\phi)$ for
 $a\in \Bbb K$. Indeed,
\bea
(w, f^*(a\phi))_W =(f(w), a\phi)^\si_V = a^\si (f(w), \phi)^\si_V
=a^\si (w, f^*(\phi))_W=(w, a^\si f^*(\phi))_W\,.
\eea

 Let  $\mathcal{C} = \{ w_1,\dots,w_n\}$ be a basis of $W$ 
and $\mathcal{B} = \{ v_1,\dots,v_m\}$ a basis of $V$. 
Let $\mc C^*=\{\phi_1, \dots, \phi_n\}$ be the dual basis of $W^*$
and  
$\mc B^*=\{\psi_1, \dots, \psi_m\}$ the dual basis of $V^*$.

If $[f]_{{\mathcal{B}} \leftarrow {\mathcal{C}}} = (a_{ij})$, then
\bean
\label{dual ma}
[f^*]_{\mc C^*\leftarrow \mc B^*} \, = \,
 ([f]_{\mc B\leftarrow\mc C}^\si) ^T\,.
\eean

\subsection{Field $\Bbb K(u)$}
In this paper we consider some particular fields $\Bbb K(u)$.

\vsk.2>
Let $u=(u_1,\dots,u_r)$ be variables. Let $\K(u)$ be the field
of rational functions
in variables 
\bean
\label{u/p}
u_i^{1/p^{s}}\,, \qquad i=1,\dots, r, \quad s\in \Z_{>0}\,,
\eean
with coefficients in $\F_p$. Thus an element of $\K(u)$
 is the ratio of two polynomials in variables $u_i^{1/p^{s}}$ with coefficients in $\F_p$.

\vsk.2>
For any
$f(u_1,\dots,u_r) \in \K(u)$, we have
\bean
\label{sqrt p}
f(u_1,\dots,u_r)^{1/p} = f(u_1^{1/p},\dots,u_r^{1/p})\,,
\eean
cf. \cite{PF}.

\vsk.2>
The field $\K(u)$ has the Frobenius automorphism
\bean
\label{Fr}
\si\,:\, \K(u)\,\to\,\K(u),\qquad f(u) \mapsto f(u)^p\,,
\eean
and its inverse
\bean
\label{i Fr}
\tau \,:\, \K(u)\,\to\,\K(u),\qquad f(u) \mapsto f(u)^{1/p}\,.
\eean

\subsection{Curve $X$}
\label{sec d.f.} 
Recall that $n=qk+1$, see Section \ref{sec Bas}.
Consider the field $\K(z)$, 
\\
 $z=(z_1,\dots,z_n)$.
 Consider the algebraic curve $X$ over $\K(z)$ defined by the affine equation 
\bean
\label{X}
y^q = F(x,z) := (x-z_1)(x-z_2)\dots(x-z_n)\,.
\eean
The curve has genus
\bean
\label{genus}
g \,: =\, k\, \frac{q(q-1)}2\,.
\eean 
The space $\Omega^1(X)$ of regular $1$-forms on $X$ is 
the direct sum
\bean
\label{Om X}
\Omega^1 (X) \,=\, \bigoplus_{a=1}^{q-1}\, \Omega^1_a (X)\,,
\eean
where  
$\dim\, \Omega^1_a (X)\,=\, ak\,,$
and  $\Omega^1_a (X)$ has  basis :
\bean
\label{basis O}
x^{i-1}\frac{dx}{y^a}\, ,   \qquad i=1,\dots, ak\,.
\eean

\subsection{Cartier operator}
\label{ssec C}

Following \cite{AH}, we introduce the {\it Cartier operator} 
\bean
\label{Cartier}
C\,: \ \Omega^1 (X) \,\to\, \Omega^1 (X)\,,
\eean
which is $\tau$-linear.
It has block structure. For $a=1,\dots,q-1$, we have
\bean
\label{Cartier2}
C \left(\Omega^1_a (X)\right)\,\subset \,\, \Omega^1_{\eta{(a)}} (X)\,,
\eean
where $\eta : \{1,\dots,q-1\} \to \{1,\dots,q-1\}$ is the division by $p$ modulo $q$ defined in \Ref{eta}.

\vsk.2>
We define the Cartier operator by the action on the basis vectors as follows.

\vsk.2>

For  $a = 1,\dots,q-1$, formula \Ref{congr} implies $q \, | \, (\eta(a)p - a)$. Hence, for $f=1$,
\dots, $ak$, we have
\bean
\label{trick}
 x^{ak - f} \frac{ dx}{y^{a}} 
&=& x^{ak - f} y^{\eta(a)p-a} \frac{ dx}{y^{a} y^{\eta(a)p-a}} 
\,=\, x^{ak - f} (y^q)^{(\eta(a)p-a)/q}\frac{ dx}
{y^{\eta(a)p}}  
\\
\notag
&=& x^{ak - f} F(x,z)^{(\eta(a)p-a)/q} \frac{ dx}{y^{\eta(a)p}}
\,=\,
 \sum_{w\geq 0} \,\,^{a} F^w_f
 (z)\, x^w \frac{ dx}{y^{\eta(a)p}}\,,
\eean
where  ${}^{a}F^w_f (z) \in \FFF_p[z]$\,.

\vsk.2>
The Cartier operator is defined by the formula
\bean
\label{def CO}
x^{ak - f}\frac{ dx}{y^{a}} \quad \mapsto \quad
 \sum_{h=1}^{\eta(a) k} \Big({}^{a}F^{(\eta(a)k-h)p+p-1}_{f} (z)\Big)^{1/p} \ x^{\eta(a)k-h}\frac{dx}{y^{\eta(a)}}\,. 
\eean

This formula has the following meaning. If $w$ in \Ref{trick} is not of the form $lp+p-1$ for some $l$, then the summand
$^aF^w_f (z)\, x^w \frac{ dx}{y^{\eta(a)p}}$ in \Ref{trick} is ignored in the definition \Ref{def CO}, and
if $w=lp+p-1$ for some $l$, then the term
$^aF^{lp+p-1}_f (z)\, x^{lp+p-1} \frac{ dx}{y^{\eta(a)p}}$ produces the summand
\bea
\Big({}^aF^{lp+p-1}_f (z)\Big)^{1/p}\, x^{l} \frac{ dx}{y^{\eta(a)}}\,
\eea
in the definition \Ref{def CO}.
It turns out that there are exactly $\eta(a) k$ such values $w=lp+p-1$ 
and that explains the upper index $\eta(a) k$ in the sum in
\Ref{def CO}.

\vsk.2>
The coefficients $ \big({}^{a}F^{(\eta(a)k-h)p+p-1}_{f} (z)\big)^{1/p}$ form the $g\times g $-matrix of the Cartier operator
with respect to the basis $x^{i-1}dx/y^a$ in \Ref{basis O}.
The matrix is  called the {\it Cartier-Manin matrix.}

\vsk.2>
Let $\Omega^1 (X)^*$ be the space dual to $\Omega^1 (X)$.  Let 
\bean
\label{d basis}
\phi ^i_a\,,\qquad a=1, \dots , q-1\,,\qquad i=1,\dots, ka\,,
\eean
denote the basis of  $\Omega^1(X)^*$  dual to the basis $(x^{i-1}dx/y^a)$ of $\Om^1(X)$.
\vsk.2>

The map 
\bea
C^* : \Omega^1 (X)^* \to \Omega^1(X)^*
\eea
 adjoint  to the Cartier operator is $\si$-linear. 
Following Serre the matrix of the map $C^*$
 is called the {\it Hasse-Witt matrix} with respect to the basis  $(\phi^i_a)$, see \cite{AH}. The entries of the Hasse-Witt matrix are
the {\it polynomials} 
\bea
{}^{a}F^{(\eta(a)k-h)p+p-1}_{f} (z)\,\in \, \F_p[z]\,,
\eea
 see \Ref{dual ma}.

\begin{rem}
The map $C^*$ is identified with the Frobenius map 
\bea
H^1(X, \mc O(X)) \to H^1(X, \mc O(X))\,,
\eea
see for example \cite{AH}.

\end{rem}

\vsk.2>
For each $a=1,\dots,q-1$,
 consider the columns of the Hasse-Witt  matrix, corresponding to the basis vectors of $\Omega^1_{\eta(a)} (X)^*$ 
and the rows corresponding to the basis vectors of $\Omega^1_{a} (X)^*$. 
The respective block of the Hasse-Witt matrix of size
$ak \times \eta(a)k$ is denoted by  ${}^{a}\mc I$. Its entries are denoted by
\bean
\label{b C-M Y}
{}^{a}\!\mathcal{I}_f^h(z)\, : =\, {}^{a} \!F^{(\eta(a)k-h)p+p-1}_f (z)\,,
\qquad  f=1,\dots,ak\,, \quad h=1,\dots,\eta(a)k\,. 
\eean

\begin{lem}
The entry ${}^{a}\mathcal{I}_f^h(z)$ is a homogeneous polynomial in $z_1,\dots,z_n$ of degree
\bea
(\eta(a)p-a)/q - (f-1) + (h-1)p\,.
\eea
\qed
\end{lem}

\subsection{Cartier map and arithmetic solutions} 
Let $W(X)$ be the $n$-dimensional $\K(z)$-vector space spanned by the  following differential 1-forms on $X$:
\bean
\label{W space}
\frac 1{x-z_i}\,\frac{dx}{y},\qquad i=1,\dots,n\,.
\eean
Notice that these are the differential 1-forms on $X$, which appear in the construction
of the solutions of the KZ equations over $\C$. Notice also that they are not regular on $X$.

\vsk.2>

Let $W(X)^*$ be the space dual to $W(X)$ and $\psi_j$, $j=1,\dots,n$,  the basis of $W(X)^*$ dual to the basis
$(dx/y(x-z_i))$.

\vsk.2>
Define  the $\tau$-linear Cartier map  
\bean
\label{main C}
\hat C\,:\, W(X)\,\to\, \Om_{a_1}^1(X)
\eean
 in the standard way. Namely we have
\bea
\frac 1{x-z_j} \,\frac{dx}y
&=&
 \frac{F(x,z)^{(a_1 p -1)/q}}{x-z_j} \frac{dx}{y^{a_1p}}
= \frac{\Phi_p(x,z,\bar M)}{x-z_j} \frac{dx}{y^{a_1p}}
\\
&=&
 P_j(z,\bar M)\,\frac{dx}{y^{a_1p}}
= \sum_i P_j^i(z,\bar M)\,x^i\,\frac{dx}{y^{a_1p}}\,,
\eea
where $\Phi_p$, $P_j$, $P_j^i$ see in \Ref{mast p}, \Ref{P(x,z,M)}, \Ref{te}.
Define $\hat  C$  by the formula
\bean
\label{def ti C}
\frac 1{x-z_j} \,\frac{ dx}{y} \quad \mapsto \quad
 \sum_{h=1}^{a_1k} \Big(P_j^{(a_1k-h)p+p-1}(z,\bar M)\Big)^{1/p} x^{a_1 k-h}\frac{dx}{y^{a_1}}\,,
\eean
cf. \Ref{def CO}.

\vsk.2>
The map 
\bean
\label{main map}
\hat  C^*\, :\, \Omega^1_{a_1} (X)^* \,\to\, W(X)^*\,,
\eean
 adjoint  to the Cartier map
$\hat  C$ is $\si$-linear. 
The matrix  of the map $\hat  C^*$
 will be
 called the {\it Hasse-Witt matrix} with respect to the 
bases  $(\phi^i_{a_1})$ and $(\psi_j)$. The entries of the Hasse-Witt matrix are
the {\it polynomials}\  $P_j^{(a_1k-h)p+p-1}(z,\bar M)$.

\vsk.2>
For any $m=1, \dots,a_1k$, we have
\bean
\label{ti C phi}
\hat  C^*(\phi^m_{a_1})\,=\,\sum_{j=1}^n P_j^{(a_1k-m)p+p-1}(z, \bar M)\,\psi_j\,.
\eean
This formula shows that
the coordinate vector 
\bean
\label{coor ve}
(P_1^{(a_1k-m)p+p-1}(z,\bar M), \dots,P_n^{(a_1k-m)p+p-1}(z,\bar M))
\eean
 of  $\hat  C^*(\phi^h_{a_1})$
is exactly the arithmetic solution $I^m(z)$ defined in \Ref{Imz}.

This construction gives us the following theorem.

\begin{thm}
\label{thm C ar1}
The Hasse-Witt matrix of the map $\hat  C^*$  defines an isomorphism 
\bean
\label{iota}
\iota_\La\, : \,\Om_{a_1}(X)^* \ \to\  \mc M_\La, \qquad \phi^m_{a_1}\ \mapsto\ I^m(z)\,,
\eean
 of the vector space $\Om_{a_1}(X)^*$ and the module $\mc M_\La$ of arithmetic solutions
for $\La=(1,\dots,1)$.

\end{thm}

\begin{proof}
Clearly the map $\iota_\La$ is an epimorphism. 
The fact that $\iota_\La$ is an isomorphism follows 
from the fact that $\dim \mc M_\La = \dim \Om_{a_1}(X)^* = a_1k$
by Lemma \ref{lem d m1}.
\end{proof}

\subsection{Arithmetic solutions from iterates}
\label{sec afi}

Consider an iterate of the Cartier map, 
\bea
C^{\circ\, b}\circ\hat  C:=C\circ C \dots \circ C \circ \hat  C \ :\  W(X)\ \to\  \Om^1(X)\,.
\eea
The image  lies in $\Om^1_{a_b}(X)$. Choose any element $\phi^{m}_{a_{b}}$ of the basis of
$\Om^1_{a_b}(X)^*$ dual to the basis $(x^{i-1}dx/ y^{a_b})$ and express the vector
\bea
(C^{\circ\, b}\circ\hat  C)^*(\phi^m_{a_b})
\eea
in terms of the basis $(\psi_j)$ of $W(X)^*$. By Theorem \ref{thm C ar1} the coordinates of that vector is an arithmetic solution,
namely, it is the solution
\bean
\label{iter sol}
\phantom{aaaa}
{}^bI^m(z):=
\sum_{m_1, \dots, m_b} \, ^{a_b}\mc I_{m_b}^{m}(z_1^{p^b}, \dots,z_n^{p^b}) 
\cdots 
\, ^{a_1}\mc I_{m_1}^{m_2}(z_1^{p}, \dots,z_n^{p})
\cdot I^{m_1}(z_1,\dots,z_n) \,.
\eean
Here $m=1,\dots, a_bk$. We will call these arithmetic solutions the 
{\it  iterated arithmetic solutions}.

\vsk.2>
We will see these solutions in Sections \ref{ssec sum} and \ref{sec sum}.

\section{Cartier map and  module $\mc M_\La$  for more general $\La$}
\label{sec C mg}

\subsection{Curve $\tilde X$}
\label{sec tX}

Recall that $n=qk+1$, see Section \ref{sec Bas}.
Let $\tilde \La = (\tilde \La_1,\dots,\tilde \La_{\tilde n})\in \Z^{\tilde n}_{>0}$ be such that
\bea
\sum_{j=1}^{\tilde n}\tilde \La_j=n\,,\qquad \tilde \La_j<q,\qquad j=1,\dots, \tilde n\,.
\eea
This means that $\tilde \La$ is a fusion of $\La=(1,\dots,1)\in \Z^n_{>0}$, see Section \ref{sec fus}.
\vsk.2>

For $ j=1,\dots,\tilde n$, $a=1,\dots,q-1$,  denote
\bean
\label{e_i}
e_j(a):= \left \lceil{\frac{\tilde \La_ja+1}q -1}\right \rceil , \qquad e(a) = \sum_{j=1}^{\tilde n} e_j(a)\,,
\eean
where $\lceil{x} \rceil$ is the smallest integer greater than $x$ or equal to $x$.

\vsk.2>
Consider the field $\K(\tilde z)$,  $\tilde z=(\tilde z_1,\dots,\tilde z_{\tilde n})$.
 Consider the algebraic curve $\tilde X$ over $\K(\tilde z)$ defined by the affine equation 
\bean
\label{X}
y^q = \tilde F(x,\tilde z) \,:=\,(x-\tilde z_1)^{\tilde \La_1}(x-\tilde z_2)^{\tilde \La_2}\dots(x-\tilde z_{\tilde n})^{\tilde \La_{\tilde n}}\,.
\eean

The space $\Omega^1(\tilde X)$ of regular $1$-forms on $\tilde X$ is 
the direct sum
\bean
\label{Om X}
\Omega^1 (\tilde X) \,=\, \bigoplus_{a=1}^{q-1}\, \Omega^1_a (\tilde X)\,,
\eean
where  $\Omega^1_a (\tilde X)$ consists of 1-forms
\bean
\label{basis O}
u(x) \,\frac{dx}{y^a}\, ,
\eean 
such that $u(x)$ is a polynomial in $x$ of degree $< ak$, and for any $j=1,\dots, \tilde n$,
the polynomial $u(z)$ has zero at $\tilde z_j$ of multiplicity at least $e_j(a)$.

This fact is checked by writing $u(x) dx/y^a$ in local coordinates on $\tilde X$ at $\infty, \tilde z_1,\dots,\tilde z_{\tilde n}$. 
Namely, at $x=\infty$ we have
\bea
x=t^{-q}, \qquad y=t^{-n}(1 + \mc O(t))\,,
\eea
where $t$ is a local coordinate.  If $d=\deg_x u(x)$, then 
$u(x)dx/y^a = -qt^{an-dq-q-1}(1 + \mc O(t))dt$. Hence
$u(x)dx/y^a$ is regular at infinity if
$an-dq-q-1\geq 0$, which is equivalent to $d < ak$.
At $x=\tilde z_j$ we have
\bea
x-\tilde z_j = t^q, \qquad y = \on{const} t^{\tilde \La_j}(1+\mc O(t))\,,
\eea
where $t$ is a local coordinate.  If $d$ is the multiplicity of $u(x)$ at $x=\tilde z_j$, then 
$u(x)dx/y^a =\on{const} t^{-a\tilde \La_j+dq+q-1}(1 + \mc O(t))dt$. Hence
$u(x)dx/y^a$ is regular  at $x=\tilde z_j$ if
$-a\tilde \La_j+dq+q-1\geq 0$, which is equivalent to $d\geq e_j(a)$.

\vsk.2>
Therefore,
\bean
\label{dim Oa}
\dim \Om^1_a(\tilde X) = ak - e(a)\,,
\eean
and  $\Om^1_a(\tilde X)$ consists of elements
\bean
\label{u(x)}
\tilde u(x) \,\frac{dx}{y^a}\, \prod_{j=1}^{\tilde n} (x-\tilde z_j)^{e_j(a)}\,,
\eean
where $\tilde u(x)$ is an arbitrary polynomial of degree less than $ak-e(a)$.

\subsection{Rank of $\mc M_{\tilde \La}$}

Consider the module $\mc M_{\tilde \La}$ of arithmetic solutions
of the KZ equations associated with $\tilde \La$.

\begin{thm}
\label{thm d e(a)}

Let $p> n=kq+1$, then the rank $d_{\tilde \La}$ of $\mc M_{\tilde \La}$ equals $a_1k-e(a_1)$.

\end{thm}

This theorem is a generalization of Lemma \ref{lem d m1}.

\begin{proof}
We have $d_\La\,=\,\big[ \sum_{j=1}^{\tilde n}\bar M_j / p \big]$
by formula
\Ref{def dL}. Here $\bar M_j$ is the minimal positive integer solution
of the congruence
\bean
\label{bMj}
M_j\equiv -\frac{\bar \La_j}{q}
\qquad 
(\on{mod} \,p)\,.
\eean
Recall the integer $a_1$ defined by
\bean
\label{a_1 n}
1\leq a_1< q\quad \on{and}\quad  q\,\vert \,(a_1p-1)\,.
\eean
Then 
\bean
\label{M0}
\tilde M_j = \tilde \La_j\,\frac{a_1p-1}q
\eean
is another solution of the congruence in \Ref{bMj}.  Denote
\bean
\label{M^0}
\tilde M = (\tilde M_1, \dots,\tilde M_{\tilde n})\,.
\eean

\begin{lem}
\label{lem M0b}
We have 
\bean
\label{M0b}
\tilde M_j = \bar M_j + e_j(a_1)p\,.
\eean

\end{lem}

\begin{proof}
Clearly, we have $\tilde M_j= \bar M_j + lp$, where
\bea
l = \left[ \frac{\tilde \La_j (a_1p-1)/q }{p} \right] = \left[ \frac{\tilde \La_j a_1}{q} - \frac{\tilde \La_j}{qp} \right].
\eea
Since $\tilde \La_j<q<p$, we have 
\bea
\frac{\tilde \La_j}{qp} < \frac{q}{qp} < \frac1q.
\eea
Furthermore, since $\tilde \La_j, a_1 < q$ and $q$ is prime, we conclude that $\tilde \La_j a_1/q$ is not an integer. These two observations imply that 
\bea
l = \left[ \frac{\tilde \La_ja_1}{q} - \frac{\tilde \La_j}{qp} \right] = \left[ \frac{\tilde \La_ja_1}{q} \right].
\eea
On the other hand, 
\bea
e_j(a_1) = \left \lceil{\frac{\tilde \La_ja_1+1}q -1}\right \rceil  = 
\left \lceil{\frac{\tilde \La_ja_1}{q}+\frac{1}{q}}\right \rceil - 1 
= \left \lceil{\frac{\tilde \La_j a_1}q }\right \rceil  -1 = \Big{[} \frac{\tilde \La_j a_1}{q} \Big{]},
\eea
where the last two equalities follow from the fact that $\tilde \La_j a_1/q$ is not an integer. The lemma is proved.
\end{proof}

To finish the proof of Theorem \ref{thm d e(a)} we observe that
\bea
d_{\tilde \La} 
&=& \left[ \sum_{j=1}^{\tilde n}\frac{\bar M_j} p \right]
=\left[ \sum_{j=1}^{\tilde n}\Big(\tilde \La_j\frac{(a_1p-1)/q} p -e_j(a_1)\Big)\right]
\\
&=& 
\left[ n\frac{(a_1p-1)/q} p\right] - e(a_1) = a_1k-e(a_1)\,,
\eea
where we use Lemmas \ref{lem M0b} and \ref{lem d m1}.
\end{proof}

\subsection{Cartier map for $\tilde X$}
Let $W(\tilde X)$ be the $\tilde n$-dimensional $\K(\tilde z)$-vector space 
spanned by the  following differential 1-forms on $\tilde X$:
\bean
\label{W space}
\frac 1{x-\tilde z_i}\,\frac{dx}{y},\qquad i=1,\dots, \tilde n\,.
\eean

Let $W(\tilde X)^*$ be the space dual to $W(\tilde X)$ and $\psi_j$, $j=1,\dots,\tilde n$, be the basis of $W(\tilde X)^*$ dual to the basis
$(dx/y(x-\tilde z_i))$.

\vsk.2>
Define  the $\tau$-linear Cartier map  $\hat C$ of  
the space $W(\tilde X)$ to the space of differential 1-forms on $\tilde X$
 in the standard way. Namely we have
\bea
\frac 1{x-\tilde z_j} \,\frac{dx}y
=
 \frac{\tilde F(x,\tilde z)^{(a_1 p -1)/q}}{x-\tilde z_j} \frac{dx}{y^{a_1p}}
= P_j(\tilde z,\tilde  M)\,\frac{dx}{y^{a_1p}}
= \sum_l P_j^l(\tilde z,\tilde  M)\,x^l\,\frac{dx}{y^{a_1p}}\,,
\eea
where  $\tilde M$, $P_j$, $P_j^l$ see in  \Ref{M^0}, \Ref{P(x,z,M)}, \Ref{te}.
Define $\hat C$ by the formula
\bean
\label{def ti C}
\frac 1{x-\tilde z_j} \,\frac{ dx}{y} \quad \mapsto \quad
 \sum_{h=1}^{a_1k} \Big(P_j^{(a_1k-h)p+p-1}(\tilde z, \tilde M)\Big)^{1/p} x^{a_1 k-h}\frac{dx}{y^{a_1}}\,,
\eean
cf. \Ref{def CO}, \Ref{def ti C}.

\begin{thm}
\label{thm im hC} For $j=1,\dots,\tilde n$, the 1-form $\hat C(dx/y(x-\tilde z_j))$ lies in $\Om^1_{a_1}(\tilde X)$,
and hence the Cartier map maps $W(\tilde X)$ to $\Om^1_{a_1}(\tilde X)$.
\end{thm}

\begin{proof}  The polynomial
\bea
u(x):=  \sum_{h=1}^{a_1k} \Big(P_j^{(a_1k-h)p+p-1}(\tilde z, \tilde M)\Big)^{1/p} x^{a_1 k-h}\frac{dx}{y^{a_1}}\,
\eea
has degree $< a_1k$. We need to check that for any $i=1,\dots,\tilde n$, the polynomial
$u(x)$ has zero at $x=\tilde z_i$ 
of multiplicity at least $e_i(a_1)$. 

Indeed, on the one hand, we have
\bea
\tilde F(x+\tilde z_i,\tilde z)^{(a_1 p -1)/q} = x^{\bar M_i + e_i(a_1)p} 
\prod_{l\ne i}(x+\tilde z_i-\tilde z_l)^{\tilde \La_i(a_1p-1)/q}\,,
\eea
by Lemma \ref{lem M0b}. Hence in the Taylor expansion
\bean
\label{tA}
\frac{\tilde F(x+\tilde z_i,\tilde z)^{(a_1 p -1)/q}}{x-\tilde z_j} =: \sum_{l} \tilde P^l_j(\tilde z) \,x^l
\eean
we have
\bean
\label{eq 0}
\tilde P^{lp+p-1}_j(\tilde z)=0\,,\qquad  l=0, \dots, e_i(a_1)-1\,.
\eean
On the other hand, we have  
\bean
\label{x+}
\frac{\tilde F(x+\tilde z_i,\tilde z)^{(a_1 p -1)/q}}{x-\tilde z_j} =
\sum_l P_j^l(\tilde z,\tilde  M)\,(x+\tilde z_i)^l\,.
\eean
By Lucas\rq{} Theorem \ref{thm L} we have
\bean
\label{Lc}
\tilde P^{lp+p-1}(\tilde z) = \sum_{h\geq 0} \binom{l+h}h \tilde z_i^{(l+h)p}P_j^{(l+h)p+p-1} (\tilde z,\tilde  M)\,,
\eean
for any $l$. Formulas \Ref{eq 0} and \Ref{Lc} show that the polynomial 
$u(x)$ has zero at $x=\tilde z_i$ of multiplicity at least $e_i(a_1)$.
The theorem is proved.
\end{proof}

\subsection{Cartier map and arithmetic solutions}
The map 
\bean
\label{i t C}
\hat C^* \,:\, \Omega^1_a (\tilde X)^* \,\to\, W(\tilde X)^*\,,
\eean
 adjoint  to the Cartier operator
$\hat C$ is $\si$-linear. 

Elements of $\Omega^1_a (\tilde X)$ have the form
\bea
u(x) \frac {dx}{y^{a_1}} = \sum_{i=0}^{a_1k-1} u^{i} x^{i} \,\frac {dx}{y^{a_1}}\,
\eea
with suitable coefficients $u^{i}$, 
see Section \ref{sec tX}. For $m=0, \dots, a_1k-1$ define an element $\phi^m\in W(\tilde X)^*$ by the formula
\bean
\label{def psi i}
\phi^m \ :\  u(x)\, \frac {dx}{y^{a_1}} \ \mapsto\  u^{m}\,.
\eean

\vsk.2>
For any $m=1, \dots,a_1k$, we have
\bean
\label{ti C phi}
\hat C^*(\phi^m)\,=\,\sum_{j=1}^{\tilde n} P_j^{mp+p-1}(\tilde z, \tilde M)\,\psi_j\,.
\eean
This formula shows that
the coordinate 
vector 
\bean
\label{ord vec}
(P_1^{mp+p-1}(\tilde z, \tilde  M), \dots,P_{\tilde n}^{mp+p-1}(\tilde z, \tilde  M))
\eean
  of  
$\hat C^*(\phi^m)$
is an arithmetic solution constructed 
in Theorem \ref{thm Fp} and all solutions constructed in Theorem \ref{thm Fp} are of this form.

\vsk.2>
This construction gives us the following theorem.

\begin{thm}
\label{thm C ar2}
The map $\hat C^*$ adjoint to  the Cartier map $\hat C : W(\tilde X) \to \Om^1_{a_1}(\tilde X)$
defines an isomorphism 
\bean
\label{iota ti}
\phantom{aaaaaaa}
\iota_{\tilde \La} \,:\, \Om_{a_1}(\tilde X)^*  \to  \mc M_{\tilde \La}\,, \quad
\phi^m \mapsto
(P_1^{mp+p-1}(\tilde z, \tilde  M), \dots,P_{\tilde n}^{mp+p-1}(\tilde z, \tilde  M))\,,
\eean
 of the vector space $\Om_{a_1}(\tilde X)^*$ and the module $\mc M_{\tilde\La}$ of arithmetic solutions
for $\tilde \La$.

\end{thm}

\begin{proof}
Clearly the map $\iota_\La$  is an epimorphism.
The fact that $\iota_\La$ is an isomorphism follows from the fact that 
$\dim \mc M_{\tilde \La} = \dim \Om_{a_1}(\tilde X)^* = a_1k-e(a_1)$
by Theorem  \ref{thm d e(a)}.
\end{proof}

\section{Hasse-Witt matrix for curve $Y$}
\label{sec HW y}

In Section \ref{sec Cartier-Manin} we introduced the curve $X$ over the field 
$\K(z)$ and determined its Hasse-Witt matrix.  In this section we consider
the same curve over a new field $\K(\la)$, where $z$ and $\la$ are related by a fractional linear transformation,
and calculate its Hasse-Witt matrix. We will use that new Hasse-Witt matrix to relate the arithmetic solutions
of the KZ equations over $\C$ and over $\F_p$.

\subsection{Curve $Y$}
\label{sec cY}

Recall that $n=qk+1$, see Section \ref{sec Bas}.
Consider the field $\K(\la)$, 
\\
 $\la=(\la_3,\dots,\la_n)$.
 Consider the algebraic curve $Y$ over $\K(\la)$ defined by the affine equation 
\bean
\label{Y}
y^q = G(x,\la) := x(x-1)(x-\la_2)\dots(x-\la_n)\,.
\eean
The space $\Omega^1(Y)$ of regular $1$-forms on $Y$ is 
the direct sum
$\Omega^1 (Y) \,=\, \bigoplus_{a=1}^{q-1}\, \Omega^1_a (Y)$\,,
where  
$\dim\, \Omega^1_a (Y)\,=\, ak\,,$
and  $\Omega^1_a (Y)$ has  basis :
\bean
\label{basis Ola}
x^{i-1}\frac{dx}{y^a}\, ,   \qquad i=1,\dots, ak\,.
\eean

\vsk.2>

We define the Cartier operator
\bean
\label{Cartier}
C\,: \ \Omega^1 (Y) \,\to\, \Omega^1 (Y)\,
\eean
in the same way as in Section \ref{ssec C}. The operator
has block structure. For $a=1,\dots,q-1$, we have
$C \left(\Omega^1_a (Y)\right)\,\subset \,\, \Omega^1_{\eta{(a)}} (Y)$.

\vsk.2>

For $f=1$,\dots, $ak$, we have
\bean
\label{tRick}
 x^{ak - f} \frac{ dx}{y^{a}} 
= x^{ak - f} G(x,\la)^{(\eta(a)p-a)/q} \frac{ dx}{y^{\eta(a)p}}
\,=\,
 \sum_{w\geq 0} \,\,^{a}G^w_f
 (\la)\, x^w \frac{ dx}{y^{\eta(a)p}}\,,
\eean
where  ${}^{a} G^w_f (\la) \in \FFF_p[\la]$\,.
The Cartier operator is defined by the formula
\bean
\label{def cO}
x^{ak - f}\frac{ dx}{y^{a}} \quad \mapsto \quad
 \sum_{h=1}^{\eta(a) k} \left({}{}^{a}G^{(\eta(a)k-h)p+p-1}_{f} (\la)\right)^{1/p} \ x^{\eta(a)k-h}
\frac{dx}{y^{\eta(a)}}\,. 
\eean

\vsk.2>
Let $\Omega^1 (Y)^*$ be the space dual to $\Omega^1 (Y)$.  Let 
\bean
\label{d basis}
\phi ^i_a\,,\qquad a=1, \dots , q-1\,,\qquad i=1,\dots, ka\,,
\eean
denote the basis of  $\Omega^1 (Y)^*$  dual to the basis $(x^{i-1}dx/y^a)$ of $\Om^1_a(Y)$.
\vsk.2>

The map 
\bea
C^* : \Omega^1 (Y)^* \to \Omega^1(Y)^*
\eea
 adjoint  to the Cartier operator is $\si$-linear. 
The matrix of the map $C^*$
 is called the {\it Hasse-Witt matrix} with respect to the basis  $(\phi^i_a)$. 
The entries of the Hasse-Witt matrix are
the {\it polynomials} ${}^{a}G^{(\eta(a)k-h)p+p-1}_{f} (\la)$, see \Ref{dual ma}.

\vsk.2>
For each $a=1,\dots,q-1$,
 consider the columns of the Hasse-Witt  matrix, corresponding to the basis vectors of $\Omega^1_{\eta(a)} (Y)^*$ 
and the rows corresponding to the basis vectors of $\Omega^1_{a} (Y)^*$. 
The respective block of the Hasse-Witt matrix of size
$ak \times \eta(a)k$ is denoted by  ${}^{a}\mc K$. Its entries are denoted by
\bean
\label{b C-M Yn}
{}^{a}\!\mathcal{K}_f^h(z)\, : =\, {}^{a} \!G^{(\eta(a)k-h)p+p-1}_f (z)\,,
\qquad  f=1,\dots,ak\,, \quad h=1,\dots,\eta(a)k\,. 
\eean

\subsection{Example}
\label{ex CM}

For $p=5$,  $ q=3$, $n=4$, we have
\bea
&&
y^3 = x(x-1)(x-\la_3)(x-\la_4)\,,
\\
&&
\Omega^1(Y) = \Omega^1_1(Y) \oplus \Omega^1_2(Y)\,=\,
 \Big\langle \frac{dx}{y}\Big\rangle\,\oplus\, \Big\langle \frac{dx}{y^2}\,,\,  x\frac{dx}{y^2}\Big\rangle\,,
\\
&& \frac{dx}{y} \ \ \mapsto \ \Big({}^{1}\mathcal{K}^1_1\Big)^{1/p} \ \frac{dx}{y^2} + 
\Big({}^{1}\mathcal{K}^2_1\Big)^{1/p} \ \frac{xdx}{y^2}\,, 
\\
&& \frac{dx}{y^2} \ \ \mapsto \ \Big({}^{2}\mathcal{K}^1_1\Big)^{1/p} \ \frac{dx}{y}\,, 
\\
&&
 \frac{xdx}{y^2} \ \mapsto \ \Big({}^{2}\mathcal{K}^1_2\Big)^{1/p} \ \frac{dx}{y}\,,
\eea
where
\bea
&&
{}^1 \mathcal{K}^1_1 (\la) = -\la_3^3 -\la_4^3 - 9 \la_3^2 \la_4
 - 9 \la_4^2 \la_3 - 9 \la_3^2 - 9 \la_4^2 -9\la_3 - 9\la_4 - 27 \la_3 \la_4 -1, 
\\
&&
{}^1 \mathcal{K}^2_1 (\la) = 3 \la_3^2 \la_4^2 ( \la_3 \la_4 + \la_3 + \la_4), 
\\
&&
{}^2 \mathcal{K}^1_1 (\la) = -\la_3 - \la_4 -1, 
\\
&&
 {}^2 \mathcal{K}^1_2 (\la) = 1.
\eea
The Hasse-Witt matrix is
\bea
\begin{pmatrix}
                            0                & {}^1 \mathcal{K}^1_1 (\la)  &  {}^1 \mathcal{K}^2_1 (\la) \\
           {}^2 \mathcal{K}^1_1 (\la)  &                 0                 &                   0               \\
           {}^2 \mathcal{K}^1_2 (\la)  &                 0                 &                   0               \\                   \end{pmatrix}.
\eea

\subsection{Homogeneous polynomials ${}^{a}\!\mathcal{J}_f^h(z)$}
\label{sec impJ}

Change variables in the polynomial ${}^{a}\mathcal{K}_f^h(\la)$,
\bea
\la_j=\frac{z_j-z_1}{z_2-z_1}, \qquad j = 3,\dots,n\,,
\eea
and multiply the result by $(z_2-z_1)^{(\eta(a)p-a)/q + (h-1)p - (f -1)}$.

\begin{lem}
\label{lem mc KJ}
The function
\bean
\label {mc JK}
{}^{a}\mathcal{J}_f^h(z)\,:=\, 
(z_2-z_1)^{(\eta(a)p-a)/q + (h-1)p - (f -1)} \cdot
{}^{a}\mathcal{K}_f^h(\la(z))
\eean
is a homogeneous polynomial  in $z_1,\dots,z_n$ of degree $(\eta(a)p-a)/q + (h-1)p - (f -1)$. 
\qed
\end{lem}

\begin{rem}
The polynomials ${}^{a}\mathcal{J}_f^h(z)$ are entries of the Hasse-Witt matrix of the curve defined by equation
\bea
y^q = x(x-(z_2-z_1))\dots(x-(z_n-z_1))\,,
\eea
and that curve is isomorphic to the curve with  equation
\bea
y^q = (x-z_1)(x-z_2)\dots(x-z_n)\,,
\eea
which is discussed in Section \ref{sec Cartier-Manin}.

\end{rem}

\subsection{Formula for ${}^{a_s}\mathcal{K}^h_f(\la)$}
\label{ssec CMfmla}

Recall the numbers $a_s$ introduced in \Ref{a_s}, the
base $p$ expansion 
\bea
\frac{p^d-1}{q} = A_0 + A_1 p + A_2 p^2 + \dots + A_{d-1} p^{d-1}
\eea
in \Ref{p exp 3}, and the relation 
\bea
A_s = \frac{a_{s+1}p-a_s}{q}\,, \quad s =0,\dots,d-1\,,
\eea
in Lemma \ref{A and a}\,. 

\vsk.2>

For  $s\in \Z_{\geq 0}\,$, $f=1,\dots,a_sk,$\   $h=1,\dots,a_{s+1}k$, define the sets
\bean
\label{sDel fh}
\qquad {}^s \Delta_f^h  = \{(\ell_3,\dots,\ell_n)\in \Z^{n-2}_{\geq 0} \ | \  0 \leq \sum_{i=3}^n \ell_i +f -1 -(h-1)p \leq A_s, \\
\phantom{aaaaaaaaaaaaaaaa} \ell_j \leq A_s \ \on{for} \ j=3,\dots,n \}. \phantom{aaaaa} \nonumber
\eean

\begin{lem}
\label{lem 6.1}
We have 
\bean 
\label{C-M fh=}
{}^{a_s}\mathcal{K}_f^h(\la) = \sum_{(\ell_3,\dots,\ell_n) \, \in\,  {}^{s} \!\Delta_f^h} 
{}^{a_s}\mathcal{K}_{f; \ell_3,\dots,\ell_n}^h (\la)\,,  
\eean
where
\bean
\label{C-M exp}
&&
^{a_s}\mathcal{K}_{f; \ell_3,\dots,\ell_n}^h (\la_3,\dots,\la_n)\,=\,(-1)^{A_s - f+1 +(h-1)p}
\\
\notag
&&
\phantom{aaaaaaaaaaa}
\times
\, \binom{A_s}{{\sum}_{i=3}^n \ell_i + f -1 - (h-1)p} 
\prod_{i=3}^n \binom{A_s}{\ell_i}\, \la_3^{\ell_3} \dots \la_n^{\ell_n}\,.
\eean
In particular all terms\  $^{a_s}\mathcal{K}_{f; \ell_3,\dots,\ell_n}^h (\la_3,\dots,\la_n)$ are nonzero.

\end{lem}

\begin{proof}
The lemma is proved by straightforward calculation similar to the proof of Theorem \ref{thm K}.
\end{proof}

 Similar formulas
can be obtained for all entries ${}^{a}\mathcal{K}^h_f(\la)$.

\section{Comparison of solutions over $\mathbb{C}$ and $\FFF_p$}
\label{sec 8}

In this section we will
\begin{enumerate}
\item
 distinguish one holomorphic solution of the KZ equations for $\La=(1,\dots,1)$, 
 \item
  expand it into the Taylor series, 
  \item
 reduce this Taylor expansion modulo $p$, 
  \item
  observe that the reduction mod $p$ of the Taylor expansion of the distinguished solution
gives  all arithmetic solutions mod $p$,  and
conversely the arithmetic  solutions together with matrix coefficients
of the Hasse-Witt matrix determine this Taylor expansion.

\end{enumerate}

\subsection{Distinguished holomorphic solution}
\label{sec dist sol}
Consider the KZ equations \Ref{KZ} for $\La=(1,\dots,1)$ over the field $\C$.
We assume that $n=qk+1$ as in Section \ref{sec Bas}.

\vsk.2>

Recall that the solutions have the form $I^{(\gamma)}(z)=(I_1(z),\dots,I_{n}(z))$, where
\bea
I_j(z) = \int_{\ga} \frac{1}{\sqrt[q]{(t-z_1)\dots (t-z_n)}}\,\frac{dt}{t-z_j}
\eea
and $\ga$ is an oriented loop on the complex algebraic curve with equation
 \bea
y^q= (x-z_1)\dots(x-z_n)\,.
\eea
Assume that $z_3,\dots,z_{n}$ are  closer to
 $z_1$ than to $z_2$:
 \bea
 \Big| \frac{z_j-z_1}{z_2-z_1}\Big|<\frac 12, \qquad
 j=3,\dots,n\,.
 \eea
Choose $\ga$ to be the circle 
 $  \big\{ t\in \C\ |\ \big|\frac{t-z_1}{z_2-z_1}\big|=\frac12\big\}$
 oriented counter-clockwise, and multiply the vector $I(z)$ by the normalization constant
 $(-1)^{1/q}/2 \pi i$.  

\vsk.2>
We will describe the normalization procedure of the solution $I(z)$ more precisely in Section \ref{sec resc}.

\vsk.2>
 
 We call the solution $I(z)$ the {\it distinguished} solution.

\subsection{Rescaling} 
\label{sec resc}

Change variables, $t=(z_2-z_1)x+z_1$, and write
\bean
\label{IL}
I(z_1,\dots,z_n)=(z_2-z_1)^{-1/q-k} L(\la_3,\dots,\la_n),
\eean
where 
\bea
\la=(\la_3,\dots,\la_n) =\Big(\frac{z_3-z_1}{z_2-z_1},\dots,\frac{z_n-z_1}{z_2-z_1}\Big),&
\eea
$L(\la) = (L_1,\dots,L_n)$,
\bean
\label{Lj}
L_j = \frac{(-1)^{1/q}}{2 \pi i} \int _{|x|=1/2} \frac{dx}{\sqrt[q]{x(x-1)(x-\la_3)\dots (x-\la_n)}}\,\frac 1{x-\la_j}\,,
\eean
and  $\la_1 = 0, \la_2 = 1$.

\vsk.2>
The integral $L(\la)$ is well-defined at $(\la_3, \dots, \la_n)=(0.\dots,0)$ and
\bea
L_j(0,\dots,0) = \frac{(-1)^{1/q}}{2 \pi i} \int _{|x|=1/2} \frac{dx}{x^k\sqrt[q]{x-1}}\,\frac 1{x-\la_j}\,.
\eea
The $q$-valued function
\bea
\frac{(-1)^{1/q}}{x^k\sqrt[q]{x-1}}
\eea
has no monodromy over  the circle $|x|=1/2$. 
To fix the value of the integrals in \Ref{Lj} we choose 
over the circle $|x|=1/2$ that branch of the 
$q$-valued function
\bea
 \frac{(-1)^{1/q}}{\sqrt[q]{x(x-1)(x-\la_3)\dots (x-\la_n)}}\,,
\eea
which is positive at $x=1/2$ and $(\la_3, \dots, \la_n)=(0.\dots,0)$.

\vsk.2>

The function $L(\la)$ is holomorphic at the point $(\la_3,\dots,\la_n)=(0,\dots,0)$. Hence
\bean
\label{L taylor}
L(\la) = \sum_{(k_3,\dots,k_n) \in \Z^{n-2}_{\geq 0}} L_{k_3,\dots,k_n} \la_3^{k_3}\dots\la_n^{k_n} 
\eean
for suitable complex numbers $ L_{k_3,\dots,k_n}$ in a neighborhood of the point $(0,\dots,0)$.

\subsection{Taylor expansion of $L(\la)$}

\begin{lem}
\label{L te exp}
We have
\bean
\label{L_k's=}
L_{k_3,\dots, k_n} 
&=&
 (-1)^{k} \binom{-1/q}{k_3+\dots+k_n+k} \prod_{i=3}^n \binom{-1/q}{k_i} 
\\
\notag
&
\times &
\Big(1, \,-q \,\Big(\sum_{i=3}^n k_i+k\Big),\,qk_3+1,\,\dots,\,qk_n+1\Big)\,,
\eean
where the integer $k$ is defined by the equation $n=qk+1$.
\end{lem}

\begin{proof} 
Indeed,
\bea
L_1(\la)
&=& \sum_{k_3,\dots,k_n=0}^{\infty} \frac{\la_3^{k_3}\dots \la_n^{k_n}}{k_3!\dots k_n!} 
\frac{\partial^{k_3+\dots+k_n} L_1}{\partial \la_3^{k_3} \dots \partial \la_n^{k_n}} (0)
\\
&
= &
\frac{(-1)^{1/q}}{2 \pi i} \sum_{k_3,\dots,k_n=0}^{\infty} \la_3^{k_3}\dots \la_n^{k_n} (-1)^{\sum_{i=3}^n k_i}  \prod_{\substack{i=3}}^n \binom{-1/q}{k_i} 
\\
&
\times
& \int_{|x|=1/2} x^{-(n-1)/q-\sum_{i=3}^n k_i -1} (x-1)^{-1/q} dx 
\\
&
= &
(-1)^{k} \sum_{k_3,\dots,k_n=0}^{\infty} \la_3^{k_3}\dots \la_n^{k_n} \binom{-1/q}{\sum_{i=3}^n k_i + k} \prod_{\substack{i=3}}^n \binom{-1/q}{k_i}.
\eea
Similarly,
\bea
L_2(\la) = (-1)^{k} \sum_{k_3,\dots,k_n=0}^{\infty} \la_3^{k_3}\dots \la_n^{k_n} \binom{-1/q}{\sum_{i=3}^n k_i + k}
 \prod_{i=3}^n \binom{-1/q}{k_i} \Big(-q \Big( \sum_{i=3}^n k_i +k \Big) \Big)
\eea
and 
\bea
L_j(\la) = (-1)^{k} \sum_{k_3,\dots,k_n=0}^{\infty} \la_3^{k_3}\dots \la_n^{k_n} \binom{-1/q}{\sum_{i=3}^n k_i + k} \prod_{\substack{i=3}}^n \binom{-1/q}{k_i} (qk_j +1 )
\eea
for $j=3,\dots, n$. The lemma is proved.
\end{proof}

\begin{cor}
Each coefficient of the series $L(\la)$ is well-defined modulo $p$. 
\end{cor}

\begin{proof}
The corollary follows from Theorem \ref{thm 1/q}.
\end{proof}

\subsection{Coefficients $L_{{k}_3,\dots,{k}_n}$ nonzero modulo $p$}

Let $({k}_3, \dots, {k}_n) \in \mathbb{Z}_{\geq0}^{n-2}$ with
\bea
{k}_i={k}_i^0+{k}_i^1 p + \dots + {k}_i^b p^b, \qquad 0\leq {k}_i^j\leq p-1, \qquad i=3,\dots,n\,,
\eea 
the base $p$ expansions. Assume that not all numbers ${k}_i^b$, $i=3,\dots,n$, are equal to zero.

Recall  the base $p$ expansion 
\bea
\frac{p^d-1}{q} = A_0 + A_1 p + A_2 p^2 + \dots + A_{d-1} p^{d-1}\,,
\qquad A_s= \frac{a_{s+1}p-a_s}{q}\,.
\eea
Extend the sequence $(A_s)$ $d$-periodically,
\bean
\label{A_s ext}
A_{s+d} := A_{s}\,.
\eean

\begin{lem}
\label{lem prod not 0}
We have
 $\prod_{i=3}^n \binom{-1/q}{{k}_i} \not\equiv 0 \; (\on{mod}\ p)$ if and only if
\bean
\label{admis1} 
{k}_i^{s} \leq A_{s}
\eean
for $i=3,\dots,n, \ s = 0,\dots, b$. 
\end{lem}
\begin{proof}
The lemma follows from Theorem \ref{thm 1/q}.
\end{proof}

\begin{lem}
\label{lem inq}
Assume that condition \Ref{admis1} holds. Then for $s=0,\dots, b$, we have
\bean
\label{inqty}
\sum_{i=3}^n {k}_i^s + a_s k \,< \,a_{s+1}kp\,.
\eean
\end{lem}

\begin{proof}

 We have 
\bea
&&
\sum_{i=3}^n {k}_i^s\, \leq\, (n-2) A_{s} \,= \,(qk-1) A_{s} \,= \,qk \ \frac{a_{s+1} p - a_{s}}{q} - A_{s} 
\\
&&
\phantom{aaaaa}
 = a_{s+1}kp - a_{s} k - A_{s} \,< \,a_{s+1}kp - a_{s} k, 
\eea
where the last inequality follows from the inequality $0< A_s$, see \Ref{A_s est}. 
\end{proof}

Following \cite{V5} define the \emph{shift coefficients} $(m_0,\dots, m_{b+1})$ as follows. Define
 $m_0=k+1$. 
For $s=0$ formula \Ref{inqty} takes the form 
\bea
\sum_{i=3}^n {k}_i^0 + m_0 -1 < a_1 kp,
\eea
since $a_0 = 1$. Hence there exists a unique integer $m_1, 1 \leq m_1 \leq a_{1} k$, such that 
\bea
0 \,\leq \,\sum_{i=3}^n {k}_i^0 + m_0 -1 - (m_1 -1) p \,< \,p\,.
\eea
For $s=1$ formula \Ref{inqty} takes the form 
\bea
\sum_{i=3}^n {k}_i^1 + a_1 k < a_2 kp.
\eea
By construction $m_1 \leq a_1 k$, hence 
\bea
\sum_{i=3}^n {k}_i^1 + m_1 -1 < a_2 kp.
\eea
Therefore, there exists a unique integer $m_2, 1 \leq m_2 \leq a_2 k$, such that 
\bea
0\, \leq \,\sum_{i=3}^n {k}_i^1 + m_1 -1 - (m_2 -1) p\, <\, p\,,
\eea
and so on.

\vsk.2>
We will obtain $m_s$ with $1 \leq m_s \leq a_s k$ for $s=1,\dots, b$.
We define $m_{b+1}$ to be the unique integer such that $1\leq m_{b+1} \leq a_{b+1}k$ and 
\bea
0\,  \leq \,\sum_{i=3}^n {k}_i^b + m_b -1 - (m_{b+1}-1)p \, < \,p\,.
\eea

We say that a tuple $({k}_3,\dots,{k}_n)$ is \emph{admissible} with respect to  $p$ if the following inequalities hold
\bean
\label{admis11} 
&&
{k}_i^{s} \leq A_{s}\,,
\qquad\qquad \qquad\qquad\qquad
 i=3,\dots,n, \ s = 0,\dots, b\,,
\\
\label{admis2}
&&
\sum_{i=3}^n {k}_i^s + m_s -1 - (m_{s+1}-1) p\, \leq \,A_s\,, \quad
\qquad s=0,\dots,b\,,
\\
\label{ad3}
&&
m_{b+1} -1 \,\leq \, A_{b+1}\,.
\eean

\begin{lem}
\label{coef not 0}
We have
\bea
\binom{-1/q}{{k}_3+\dots+{k}_n+k} \,\not\equiv\, 0 \qquad  (\on{mod}\ p)\,,
\eea 
if and only if the shift coefficients of the tuple $({k}_3,\dots,{k}_n)$ satisfy \Ref{admis2} and \Ref{ad3}.
\end{lem}

\begin{proof}
The $p$-ary expansion of ${k}_3+\dots+{k}_n+k$ is 
\bea
&&
\left( \sum_{i=3}^n {k}_i^0 + m_0 -1 - (m_1-1)p \right) + \left( \sum_{i=3}^n {k}_i^1 + m_1 -1 - (m_2-1)p \right) p + \dots 
\\
&&
\phantom{aaa}
 + \left( \sum_{i=3}^n {k}_i^{b} + m_b -1 - (m_{b+1}-1)p \right) p^b + (m_{b+1}-1) p^{b+1}. \phantom{aaaa}
\eea
By Theorem \ref{thm 1/q} the following congruence holds modulo $p$: 
\bea
 \binom{-1/q}{{k}_3+\dots+{k}_n+k} \
 \equiv \ \binom{A_{b+1}}{m_{b+1}-1} \cdot \prod_{s=0}^{b} \binom{A_{s}}{\sum_{i=3}^n {k}_i^s + m_s -1 - (m_{s+1}-1)p}.
\eea
The right-hand side of the congruence above is nonzero,
 if and only if \Ref{admis2} and \Ref{ad3} hold.
\end{proof}

Recall  the sets ${}^{s} \Delta_{f}^{h}$  defined in \Ref{sDel fh}.

\begin{lem}
The tuple $({k}_3,\dots,{k}_n)$ is admissible if and only if 
$m_{b+1} -1 \leq A_{b+1}$ and 
\bea
({k}_3^s,\dots,{k}_n^s) \in {}^{s} \Delta_{m_s}^{m_{s+1}}
\qquad\on{for} \ s = 0,\dots, b\,.
\eea
\qed
\end{lem}

\begin{thm} 
\label{thm C-M}

The following statements hold true:
\begin{itemize}
\item[(i)] $L_{{k}_3,\dots,{k}_n} \ \not\equiv 0 \ (\on{mod}\ p)$ if and only if the tuple $({k}_3,\dots,{k}_n)$ is admissible. 
\item[(ii)] If the tuple $({k}_3,\dots,{k}_n)$ is admissible, then 
\bean
\label{L_k's C-M}
&&
 L_{{k}_3,\dots, {k}_n} \la_3^{{k}_3}\dots \la_n^{{k}_n}
 \equiv  \ (-1)^{(a_{b+1}p^{{b}+1}-1)/q+m_{b+1}-1} \binom{A_{b+1}}{m_{b+1}-1}  
 \\
&&
\phantom{aaaaa}
\times \left( \prod_{s=1}^{b} \,
{}^{a_{s}}\!\mathcal{K}_{m_s; {k}_3^s,\dots,{k}_n^s}^{m_{s+1}}
 \left( \la_3^{p^s},\dots,\la_n^{p^s} \right) \right) K^{m_1}_{{k}_3^0,\dots,{k}_n^0} (\la)  \qquad (\on{mod}\ p)\,,
\notag
\eean
where ${}^{a_{s}}\!\mathcal{K}_{f; \vec k}^{h} (\la)$ 
are terms of the Hasse-Witt matrix expansion in \Ref{C-M fh=} 
and $K^{m}_{\vec{k}}(\la)$ are the terms of the expansion of vector polynomial $K^m(\la)$  in \Ref{K=}.   
\end{itemize}
\end{thm}

\begin{proof}
We have $L_{k_3,\dots,k_n} \ \not\equiv 0 \ (\on{mod}\ p)$ if
 and only if each of the binomial coefficients in \Ref{L_k's=} is
 not divisibly by $p$. By Lemmas \ref{lem prod not 0} and 
\ref{coef not 0} this is equivalent to saying that the tuple $(k_3,\dots,k_n)$ is admissible. This gives part (i). 

By Lemma \ref{L te exp} we have
\bean
\label{again L}
 L_{{k}_3,\dots, {k}_n} 
 &=& (-1)^{k} \binom{-1/q}{{k}_3+\dots+{k}_n+k} \prod_{i=3}^n \binom{-1/q}{{k}_i} 
\\
&\times &
\Big{(}1,\, -q \,\Big(\sum_{i=3}^n {k}_i+k\Big),\,q{k}_3+1,\,\dots,\,q{k}_n+1\Big{)}
 \,. \nonumber
\eean
Theorem \ref{thm 1/q} allows us to write the binomial coefficients in 
formula \Ref{again L} as products and then formula \Ref{L_k's C-M} 
becomes a straightforward corollary of formulas for ${}^{a_{s}}\mathcal{K}^h_{f;\vec{k}}$ and $K^m_{\vec{k}}$.

Notice that calculating the power of -1 in the right-hand side of formula \Ref{L_k's C-M}
we use the identity
\bea
A_0+A_1p+\dots+A_bp^b =
\frac{a_1p-1}q + \frac{a_2p-a_1}qp+\dots + \frac{a_{b+1}p-a_b}qp^b = \frac{a_{b+1}p^{b+1}-1}q\,.
\eea
\end{proof}

\subsection{Decomposition of $L(\la)$ as a sum of $K^m(\la)$}
\label{ssec sum}

Define the set 
\bean
\label{set M}
M = \{ (m_0,\dots,m_{b+1}) \in \mathbb{Z}^{b+2}_{\geq 1} \ | \ b \in \mathbb{Z}_{\geq 0},  \ m_0=k+1, \ 1\leq m_s \leq a_{s} k \\
\phantom{aaa} \on{for} \ s = 1,\dots, b, \ 1\leq m_{b+1} \leq a_{b+1} k \}. \nonumber
\eean
For any $\vec{m}=(m_0,\dots,m_{b+1}) \in M$, define the $n$-vector of polynomials in $\la$\,: 
\bean
\label{N}
N_{\vec{m}} (\la) 
&=& (-1)^{(a_{b+1}p^{{b}+1}-1)/q+m_{b+1}-1} \binom{A_{b+1}}{m_{b+1}-1} 
\\
&
\times&
 \left( \prod_{s=1}^b\  {}^{a_{s}}\! \mathcal{K}^{m_{s+1}}_{m_s} \left( \la_3^{p^s},\dots,\la_n^{p^s} \right) \right)  K^{m_1} (\la_3,\dots,\la_n). \nonumber
\eean

\begin{thm}
\label{thm dec}
We have
\bean
\label{dec la}
L(\la) \ \equiv\  \sum_{\vec m \in M} N_{\vec m}(\la) \qquad (\on{mod}  p)\,.
\eean
Moreover, if a monomial $\la_3^{k_3}\dots\la_{n}^{k_{n}}$ enters  one of the vector polynomials
$N_{\vec m}(\la)$ with a 
 nonzero coefficient,  then this monomial does not enter with nonzero coefficient any other vector polynomial
$N_{\vec m\rq{}}(\la)$.

\end{thm}

\begin{proof}
The theorem is a straightforward corollary of Theorem \ref{thm K}, Lemma \ref{lem 6.1}, and Theorem \ref{thm C-M}.
\end{proof}

For $q=2$ this theorem is \cite[Corollary 7.4]{V5}.

\subsection{Distinguished solution over $\C$ and solutions $J^m(z)$ over $\F_p$} 
\label{sec sum}

Consider the distinguished  solution $I(z_1,\dots,z_n)$
of the KZ equations \Ref{KZ} over 
$\C$, see \Ref{IL},  and the arithmetic solutions $J^m(z)$ of the same equations over $\F_p$, see \Ref{J&I}.

We have
\bean
\label{IL n}
I(z_1,\dots,z_n)
&=&(z_2-z_1)^{-1/q-k} \cdot L(\la_3,\dots,\la_n)\,,
\\
\label{J&K n}
J^m(z) 
&=&
 (z_2-z_1)^{(a_1p-1)/q + (m-1)p-k}\cdot K^m(\la)\,,
\\
\label {mc JK n}
{}^{a}\mathcal{J}_f^h(z) 
&=& (z_2-z_1)^{(\eta(a)p-a)/q + (h-1)p - (f -1)} \cdot
{}^{a}\mathcal{K}_f^h(\la(z))\,,
\eean
see \Ref{IL}, \Ref{J&K}, and \Ref{mc JK}.

\vsk.2>
Expressing  $L(\la_3,\dots,\la_n)$, 
$ K^m(\la)$,  $\mathcal{K}_f^h(\la(z))$ in terms of 
$I(z_1,\dots,z_n)$, $J^m(z) $,  ${}^{a}\mathcal{J}_f^h(z)$\,
from these equations, and using the congruence
\bean
\label{dec la n}
L(\la) \ \equiv\  \sum_{\vec m \in M} N_{\vec m}(\la) \qquad (\on{mod}  p)\,
\eean
of Theorem \ref{thm dec}, we obtain a relation between the 
distinguished holomorphic solution $I(z)$ and the arithmetic data
$J^m(z) $ and ${}^{a}\mathcal{J}_f^h(z)$.

This relation shows that knowing the distinguished holomorphic solution we may recover
the arithmetic solutions as well as the entries of the associated Hasse-Witt matrix.
Conversely, knowing the arithmetic solutions and the entries of the associated Hasse-Witt matrix, we may 
recover the distinguished holomorphic solution modulo $p$. 

\begin{rem}
Notice that summands in the right-hand side of \Ref{dec la} correspond to the iterated arithmetic 
solutions defined in Section \ref{sec afi}. We may interpret formula \Ref{dec la} as a statement that
the Taylor expansion of the distinguished holomorphic solution
reduced modulo $p$ is the sum of all iterated arithmetic solutions.

\end{rem}

\end{document}